\documentclass[12pt]{article}
\usepackage[margin=2cm]{geometry}
\usepackage{amsfonts}
\usepackage{amssymb}
\usepackage{amsthm}
\usepackage{amsmath}
\usepackage{latexsym}
\usepackage{color}
\usepackage{mathrsfs}

\begin{document}
\allowdisplaybreaks[4]
\newtheorem{lemma}{Lemma}
\newtheorem{pron}{Proposition}
\newtheorem{thm}{Theorem}
\newtheorem{Corol}{Corollary}
\newtheorem{exam}{Example}
\newtheorem{defin}{Definition}
\newtheorem{remark}{Remark}
\newtheorem{property}{Property}
\newcommand{\la}{\frac{1}{\lambda}}
\newcommand{\sectemul}{\arabic{section}}
\renewcommand{\theequation}{\sectemul.\arabic{equation}}
\renewcommand{\thepron}{\sectemul.\arabic{pron}}
\renewcommand{\thelemma}{\sectemul.\arabic{lemma}}
\renewcommand{\thethm}{\sectemul.\arabic{thm}}
\renewcommand{\theCorol}{\sectemul.\arabic{Corol}}
\renewcommand{\theexam}{\sectemul.\arabic{exam}}
\renewcommand{\thedefin}{\sectemul.\arabic{defin}}
\renewcommand{\theremark}{\sectemul.\arabic{remark}}
\renewcommand{\theproperty}{\sectemul.\arabic{property}}
\def\REF#1{\par\hangindent\parindent\indent\llap{#1\enspace}\ignorespaces}

\title{\large\bf Precise large deviations of some risk objectives related to the net loss process in two nonstandard risk models}
\author{\small Yang Chen$^1$\thanks{Research supported by Humanities and Social Sciences Foundation of the Ministry of Education of China (No. 18YJC910004).},~~Zhaolei Cui$^2$~~and~~Yuebao Wang$^3$
\thanks{Corresponding author.
Telephone: +86 512 67422726. Fax: +86 512 65112637. E-mail:
ybwang@suda.edu.cn}%~~and~~Kaiyong Wang$^3$
\\
{\footnotesize\it 1 School of Mathematical Sciences, Suzhou University of Science and Technology, Suzhou, 215009, China}\\
{\footnotesize\it 2 School of mathematics and statistics, Changshu Institute of Technology, Suzhou 215000, China}\\
{\footnotesize\it 3 School of Mathematical Sciences, Soochow University, Suzhou 215006, China}\\
%{\footnotesize\it 3 School of mathematics Sciences, Suzhou University of science and technology, Suzhou 215009, China}
}
\date{}

\maketitle {\noindent\small {\bf Abstract }}\\

For two nonstandard renewal risk models,  we investigate the precise large deviations of the finite-time ruin probability
and a random sum of the net-loss process, and the asymptotics of the random-time ruin probability. Notably, in one of these models,  claim sizes series and claim interval time series are allowed to be arbitrarily dependent.
Subsequently, we apply these results to obtain precise large deviations of proportional-net-loss process
and excess-of-net-loss process,
as well as asymptotic estimates of the mean of stop-net-loss reinsurance treaty.
These results all involve the income items of the risk model, which are relatively rare in the existing references.\\



\section{\normalsize\bf Preliminary}
\setcounter{equation}{0}\setcounter{thm}{0}\setcounter{lemma}{0}\setcounter{remark}{0}\setcounter{pron}{0}\setcounter{Corol}{0}           %(Section 1)

It is widely recognized that in a risk model,  the initial capital and the duration of operation of an insurance or reinsurance business are significant. Therefore, it is important to estimate the asymptotic properties of certain risk objects when both prerequisites tend to infinity together.
As a result, precise large deviation naturally becomes an important focus of this research.
For a brief review of relevant research, please refer to Section 2.
In this section, given the close relationship between the research and the properties of distributions of relevant random variables and their dependent structure, we will introduce the concepts and marks of related distribution class, dependent structure of random variables, and two nonstandard renewal risk models before presenting the main results of this paper.

In this paper, all limit relations refer to $x\to\infty$
and all distributions are supported on $(-\infty,\infty)$ or $[0,\infty)$ without a special statement.
In addition, the following marks and conventions are used throughout the text.

Let two functions $g_1(\cdot)$ and $g_2(\cdot)$ on $[0,\infty)$ be positive eventually.
For $1\le i\neq j\le2$, we set
$$g_{i,j}=\limsup g_i(x)g_j^{-1}(x).$$
Then $g_i(x)=O\big(g_j(x)\big)$ means $g_{i,j}<\infty$, $g_1(x)\asymp g_2(x)$ means $\max\{g_{1,2},g_{2,1}\}<\infty$,
$g_i(x)\lesssim g_j(x)$ means $g_{i,j}\le1$, $g_1(x)\sim g_2(x)$ means $g_{1,2}=g_{2,1}=1$, and $g_i(x)=o\big(g_j(x)\big)$ means $g_{i,j}=0$.

Let $F$ be a distribution. We denote the tail of $F$ by $\overline{F}=1-F$,
and the $n$-fold convolution of $F$ with itself by $F^{*n}$ for all $n\ge2$.

Other marks and conventions will be given in succession below.

\subsection{\normalsize\bf Some distribution classes}

We say that the distribution class
\begin{eqnarray*}
{\cal L}=\big\{F:\overline{F}(x-t)\sim\overline{F}(x)\ \text{for each}\ t\in(-\infty,\infty)\big\}
\end{eqnarray*}
is long-tailed, and the distribution class
\begin{eqnarray*}
{\cal S}=\big\{F:F\in\mathcal{L}\ \ \text{and}\ \ \overline{F^{*2}}(x)\sim2\overline{F}(x)\big\}.
\end{eqnarray*}
is subexponential, which was introduced by Chistyakov \cite{C1964}.
Further, Lemma 2 of this paper notes that the requirement $F\in\mathcal{L}$ is unnecessary if $F$ is supported on $[0,\infty)$.
For a detailed discussion of classes $\mathcal{S}$ and $\mathcal{L}$, see Embrechts et al. \cite{EKM1997},
Resnick \cite{R2007}, Borovkov and Borovkov \cite{BB2008}, Foss et al. \cite{FKZ2013}, etc.
In particular, if $F\in\mathcal{L}$, then the following set of positive functions is not empty:
\begin{eqnarray}\label{100}
\mathcal{H}_F=\Big\{h(\cdot)\ \text{on}\ [0,\infty):h(x)\uparrow\infty,h(x)x^{-1}\downarrow0
\ \text{and}\ \int_{h(x)}^{x-h(x)}\overline{F}(x-y)F(dy)=o\big(\overline{F}(x)\big)\Big\}.
\end{eqnarray}

In addition, for each $y\in(0,\infty)$, we set
$$\overline{F_{*}}(y)=\liminf{\overline F}(xy){\overline F}^{-1}(x),
~~\overline{F^{*}}(y)=\limsup{\overline F}(xy){\overline F}^{-1}(x)\ \ \text{and}\ \ L_F=\lim\limits_{y\downarrow1}\overline{F_*}(y).$$
Then we say that the distribution classes
\begin{eqnarray*}
\mathcal{D}=\big\{F:{\overline{F_*}}(y)>0\ \ \text{for each}\ y\in(1,\infty)\big\},\ \ {\cal C}=\{F:L_F=1\},
\end{eqnarray*}
\begin{eqnarray*}
\mathcal{ERV}(\alpha,\beta)=\{F:y^{-\beta}\leq\overline{F_*}(y)\leq\overline{F^*}(y)\leq y^{-\alpha}
\ \ \text{for each}\ y>1\}
\end{eqnarray*}
for each pair $0<\alpha\le\beta<\infty$ and
$$\mathcal{ERV}=\bigcup_{0<\alpha\le\beta<\infty}\mathcal{ERV}(\alpha,\beta)$$
are dominated varying tailed, consistently varying tailed,
extended regularly varying tailed with indexes $0<\alpha\le\beta<\infty$ and extended regularly varying tailed, respectively.
Particularly, if $\alpha=\beta$, then $\mathcal{ERV}(\alpha,\beta)$ reduces to the regularly varying tailed distribution class, denoted
by $\mathcal{R}_{\alpha}$.
Some properties of these distributions are introduced by the following two propositions.

\begin{pron}\label{Proposition101}
%$(i)$ The following relations hold:\\
%$$F\in\mathcal{D}\Longleftrightarrow\overline{F}^*(y)<\infty\ \text{for each}\ y\in(0,1),$$
%where the mark $``\Longleftrightarrow"$ means ``if and only if".

%$(ii)$
The following inclusion relations are proper: for each pair $0<\alpha\le\beta<\infty$,
$$\mathcal {R}_{\alpha}\subset \mathcal{ERV}(\alpha,\beta)\subset\mathcal{C}\subset\mathcal{L}\cap\mathcal{D}\subset\mathcal{S}\subset\mathcal{L}.$$
\end{pron}

Further, we denote the moment index of $F$ by
$$I_F=\sup\Big\{s:\int_0^\infty y^sF(dy)<\infty\Big\},$$
and the upper Matuszewska index and lower Matuszewska index of distribution $F$ by
\begin{eqnarray*}
J_F^+=-\lim\limits_{y\to\infty}\ln \overline{F_*}(y)\ln^{-1} y\ \ \ \ \text{and}\ \ \ \
J_F^-=-\lim\limits_{y\to\infty}\ln \overline{F^*}(y)\ln^{-1} y.
\end{eqnarray*}

\begin{pron}\label{Proposition102} Let $F$ be a distribution in class $\mathcal{D}$.

$(i)$ $0\le J_F^-\le I_F\le J_F^+<\infty$.

$(ii)$ For each $p<J_F^-$, $\overline{F}(x)=o(x^{-p})$.

$(iii)$ For each $p>J_F^+$, $x^{-p}=o\big(\overline{F}(x)\big)$.
%and there exist two positive constants $C_1=C_1(F,p)$ and $C_2=C_2(F,p)$ such that \begin{eqnarray*} \overline{F}(x)\overline{F}^{-1}(y)\le C_2x^py^{-p}\ \ \ \ \ \ \text{for each pair}\ \ y\ge x\ge C_1. \end{eqnarray*}
\end{pron}

The above concepts and properties can be found in some references,
such as Feller \cite{F1971}, Bingham et al. \cite{BGT1987}, Cline and Samorodnitsky \cite{CS1994}, Embrechts et al. $\cite{EKM1997}$,
Kl\"{u}ppelberg and Mikosch $\cite{KM1997}$ and Tang and Tsitsiashvili $\cite{TT2003}$.

\subsection{\normalsize\bf Some dependent structures}

Wang et al. $\cite{WWG2013}$ introduced the concept of widely dependence structure of random variables.
By definition, $X_i,i\ge 1$ are said to be widely upper orthant dependent (WUOD),
if for each $n\ge1$, there exists some positive number $g_{U,F}(n)$ such that,
\begin{eqnarray}\label{101}
P\Big(\bigcap^{n}_{i=1}\{X_i > x_i\}\Big)\leq g_{U,F}(n)\prod_{i=1}^nP(X_i> x_i),\ \ \ \ \ \ x_i\in(-\infty,\infty),\ \ \ 1\le i\le n;
\end{eqnarray}
they are said to be widely lower orthant dependent (WLOD),
if for each $n\ge1$, there exists some  positive number $g_{L,F}(n)$ such that,
\begin{eqnarray}\label{102}
P\Big(\bigcap^{n}_{i=1}\{X_i \leq x_i\}\Big)\leq g_{L,F}(n)\prod_{i=1}^n P(X_i\leq x_i),\ \ \ \ \ x_i\in(-\infty,\infty),\ \ \ 1\le i\le n;               %(1.3)
\end{eqnarray}
and they are said to be widely orthant dependent (WOD) if they are both WUOD and WLOD.

WUOD, WLOD and WOD structures can be called  widely dependent (WD) as a joint name.
And $g_{U,F}(n),~g_{L,F}(n),\ n\geq 1$, are called dominating coefficients.
Clearly,
$$g_{U,F}(n)\ge1,\ \ \ g_{L,F}(n)\ge1,\ \ n\ge2\ \ \ \text{and}\ \ \ g_{U,F}(1)=g_{L,F}(1)=1.$$

Some basic properties of WD random variables are as follows, see Proposition 1.1 of Wang et al. $\cite{WWG2013}$.

\begin{pron}\label{Proposition103}
$(i)$ Let $X_i,i\geq1$ be WLOD (or WUOD). If $f_i(\cdot),i\geq1$ are nondecreasing,
then $f_i(X_i),i\geq1$ are still WLOD (or WUOD);
if $f_i(\cdot),i\geq1$ are nonincreasing, then $f_i(X_i),i\geq1$ are WUOD (or WLOD).

$(ii)$ If $X_i,i\geq1$ are nonnegative and WUOD, then for each $n\geq1$,
$$E\prod^{n}_{i=1}X_{i} \leq g_{U,F}(n)\prod^{n}_{i=1}EX_{i}.$$
In particular, if $X_i,i\geq1$ are WUOD, then for each $n\geq1$ and each $s>0$,
$$Ee^{s \sum^{n}_{i=1}X_{i}} \leq g_{U,F}(n)\prod^{n}_{i=1}Ee^{s {X_i}}.$$
\end{pron}

Further, Wang et al. $\cite{WWG2013}$ provided some examples of WD random variables,
which showed that the WD structure may can include common negatively dependent random variables,
some positively dependent random variables, and some others.

When $g_{U,F}(n)=g_{L,F}(n)=M$ for all $n\ge1$ and some $M>0$,
the inequalities (\ref{101}) and (\ref{102}) describe extended negatively upper
and lower orthant dependent (ENUOD and ENLOD) random variables, respectively.
$X_i,i\ge 1$ are said to be extended negatively orthant dependent (ENOD) if they are both ENUOD and ENLOD.
ENOD, ENUOD, ENLOD random variables are collectively called END r.v.s, see Liu $\cite{L2009}$.

Further, if $M=1$, then we have the corresponding notions of NUOD, NLOD, NOD and ND random variables,
see, for example, Ebrahimi and Ghosh $\cite{EG1981}$ and Block et al. $\cite{BSS1982}$.

In particular, we say that random variables $X_i,i\ge1$ are negatively associated (NA),
if for any disjoint nonempty subsets $A$ and
$B$ of $\{1,\cdots,m\},\ m\geq2$ and any coordinate-wise
nondecreasing functions $f_1(\cdot)$ and $f_2(\cdot)$, the inequality
\begin{eqnarray*}
Cov(f_1(X_i:i\in A),f_2(X_j:j\in B))\leq0
\end{eqnarray*}
holds whenever the moment involved exists.
For further details, please refer to Joag-Dev and Proschan \cite{JP1983}.

%In addition to Wang et al. $\cite{WWG2013}$, many literatures give the research results on WD random variables,
%see Wang and Cheng $\cite{WC2011,WC2020}$ for renewal theory,
%Wang et al. $\cite{WYL2012}$ and He et al. \cite{HCW2013} for precise large deviation theory,
%Qiu and Chen $\cite{QC2014}$, Wang et al. $\cite{WXHV2014}$, Naderi et al. \cite{NMAB2015}, Wang and Hu $\cite{WH2015}$, Chen et al. \cite{CWC2016},
%Istv\'{a}n et al. \cite{IPB2018}, Deng and Wang \cite{DW2018}, Chen and Sung \cite{CS2019}, Wu et al. \cite{WWR2019}, Sen et al. \cite{STW2020},
%Xi et al. \cite{XWCW2020} and Deng and Wang \cite{DW2020} for limit theory and statistical theory, among others.

%Besides, there are also many results on ND and END random variables, which will not be described in detail here.

\subsection{\normalsize\bf Nonstandard renewal risk models}

Here, we introduce some concepts and properties of renewal risk model for the insurance business.

In this model, let claim sizes $Y_i,\ i\geq1$ be random variables
with common distribution $G$ on $[0,\infty)$ and mean $0<\mu_G<\infty$;
the inter-arrival times $Z_i,\ i\geq1$ of the claim sizes also be random variables
with common distribution $H$ on $[0,\infty)$ and mean $0<\mu_H<\infty$.
In addition, let random time $\tau$ be a nonnegative random variable.

The times of successive claims, $S^H_n=\sum_{i=1}^nZ_i,\ n\geq1$, constitute a counting process
$$\Big(N(t)=\sum_{n=1}^\infty\textbf{1}_{\{S^H_n\leq\  t\}}:t\ge0\Big)$$
with the corresponding mean function $\lambda(t)=EN(t)$.
For convenience, we might as well assume that
$$\inf\{t:0<\lambda(t)<\infty\}=0.$$
In addition, let $c$ be a positive constant interest rate and $x$ be a nonnegative initial capital of the insurance business.
In order to ensure normal operation of the insurance business, a safe load condition is usually assumed that
$$\mu_G<c\mu_H.$$

The aforementioned model is referred to as the standard renewal risk model.
if the random vectors $(Y_i,Z_i),i\geq1$ are independent and identically distributed,
and $\{Y_i,\ i\geq1\}$, $\{Z_i,\ i\geq1\}$ and $\tau$ are mutually independent.
Under these assumptions, we define $\big(N(t):t\ge0\big)$ as the standard renewal counting process and $\lambda(\cdot)=EN(\cdot)$ as the standard renewal function.
Alternatively, if at least one of these independent assumptions is not true,
the model is considered nonstandard and $\big(N(t):t\ge0\big)$.
In this case, $\lambda(\cdot)$ are called quasi-renewal counting process and quasi-renewal function, respectively.
They are collectively called the renewal risk model, the renewal counting process and the renewal function, respectively.

In a renewal risk model, we call $X_i=Y_i-cZ_i$ the net-loss when the $i$th claim comes with distribution $F$
on $(-\infty,\infty)$ for $i\ge1$.
Further, we call the process
$$(R_{00}(t)=S^G_{N(t)}=\sum_{i=1}^{N(t)}Y_i:t\ge0)$$
the total claim amount at time $0\le t<\infty$, see Subsection 1.3.3 of Embrechts et al. \cite{EKM1997},
and the process
$$\Big(R_0(t)=\sum_{i=1}^{N(t)}(Y_i-cZ_i)=\sum_{i=1}^{N(t)}X_i=S^F_{N(t)}=S^G_{N(t)}-cS^H_{N(t)}:t\ge0\Big)$$
a net loss process. Then the maximum of sums of net loss at time $t$ constitutes a random process
$$\Big(\widetilde{R_0}(t)=\max\limits_{0\leq n\leq N(t)}\sum_{i=1}^n(Y_i-cZ_i):t\ge0\Big).$$
The process $(R_{00}(t):t\ge0)$ is the basic research object in the existing references,
and the processes $(\widetilde{R_{0}}(t):t\ge0)$ and $(R_{0}(t):t\ge0)$ are the two primary research objects of this paper.

Naturally, the finite-time ruin probability at time $t\ge0$ for some insurance business is defined by
\begin{eqnarray}\label{403}
\psi(x;t)=P\big(\widetilde{R_0}(t)>x\big),\ \ \ \ \ x\ge0,
\end{eqnarray}
Therefore, the random time ruin probability is defined by
\begin{eqnarray}\label{404}
\psi(x;\tau)=P\big(\widetilde{R_{0}}(\tau)>x\big),\ \ \ \ x\ge0.
\end{eqnarray}

Here, we can appreciate the importance of studying the net-loss process, as it constitutes a key component in determining the above-mentioned two ruin probabilities, as well as other risk metrics, see, for example, Section 4.

For the renewal risk model, some renewal theorems play important roles.
The following three results for the quasi-renewal counting process $N(t)$
generated by WD random variables $Z_i,i\ge1$ are presented as below.
The first result comes from Lemma 2.2 of Wang and Cheng \cite{WC2011}.

\begin{pron}\label{pron104}
In a nonstandard renewal risk model,
let $Z_i,\ i\ge1$ be WLOD random variables with dominating coefficients $g_{L,H}(n),\ n\ge1$.
For some $a>0$, if
\begin{eqnarray}\label{103}
\lim_{n\to\infty}g_{L,H}(n)e^{-an}=0,
\end{eqnarray}
then for each $\delta>0$, there exists $r=r(H,a,\delta)>0$ such that
\begin{eqnarray}\label{104}
\lim_{t\to\infty}Ee^{rN(t)}\textbf{\emph{1}}_{\{N(t)>(1+\delta)\mu_H^{-1}t\}}=0.
\end{eqnarray}
\end{pron}

The second result dues to Lemma 4.2 of Wang et al. \cite{WCWM2012}.

\begin{pron}\label{pron105}
In a nonstandard renewal risk model, let $Z_i,\ i\ge1$ be WLOD random variables with the dominating coefficients $g_{L,H}(n),\ n\ge1$.
For some $b>0$, if
\begin{eqnarray}\label{105}
\lim_{n\to\infty}g_{L,H}(n)n^{-b}=0,
\end{eqnarray}
then for each $k\ge1$,
\begin{eqnarray}\label{1060}
\lim_{t\to\infty}EN^k(t)(\mu_H^{-1} t)^{-k}=1.
\end{eqnarray}
\end{pron}

The final result belongs to Theorem 2.4 of Wang and Cheng \cite{WC2020} or Theorem 4 of Chen et al. \cite{CWC2016},
which makes up a gap in the proof of Theorem 1.4 of Wang and Cheng \cite{WC2011}.
In the following, we say that a positive function $g(\cdot)$ is almost decreasing (or nonincreasing),
denoted by $g(x)\ \overline{\downarrow}$,
if there exist two positive finite constants $x_0=x_0\big(g(\cdot)\big)$ and $C=C(x_0)$ such that
$$\sup_{x_0\le x\le y}g(y)g^{-1}(x)\le C.$$
Similarly, we say that a positive function $g(\cdot)$ is almost increasing (or nondecreasing),
denoted by $g(x)\ \underline{\uparrow}$, if $g^{-1}(\cdot)$  is almost decreasing.

\begin{pron}\label{pron106}
In a nonstandard renewal risk model, let $Z_i,\ i\ge1$ be WOD random variables
with dominating coefficients $g_{U,H}(n)$ and $g_{L,H}(n),\ n\ge1$.
If there exists some almost increasing function $g(\cdot)$,
some $l\ge1$ and $\omega\in(0,1)$ such that $EZ_1g(Z_1)<\infty$,
\begin{eqnarray}\label{107}
\frac{x^{l-1}}{g(x)}\ \overline{\downarrow}\ 0,\ \ \ \frac{g(x)}{x^{l-1+\omega}}\ \overline{\downarrow}\ 0
\ \ \ \text{and}\ \ \ \max\big\{g_{U,H}(n),g_{L,H}(n)\big\}\le g(n),\ \ \ n\ge1,
\end{eqnarray}
then
\begin{eqnarray}\label{108}
&\lim_{n\to\infty}n^{-1}S_n^H=\mu_H\ \ \ \ \text{and}\ \ \ \ \lim_{t\to\infty}N(t)\mu_Ht^{-1}=1,\ \ \ \ a.s.
\end{eqnarray}
%and for each $p>0$ and any $\delta>0$,
%\begin{eqnarray}\label{3020}
%EN^p(t)\textbf{\emph{1}}_{\{N(t)>(1+\delta)\lambda(t)\}}=o\big(\lambda(t)\big).
%\end{eqnarray}
\end{pron}

\begin{remark}\label{rem101}
$(i)$ In Proposition \ref{pron106}, if $Z_i,\ i\ge1$ are ENOD random variables,
then we only need the condition $EZ_1<\infty$.
At this time, there is a positive function $g_0(\cdot)$ such that $g_0(x)\underline{\uparrow}\infty$,
$EZ_1g_0(Z_1)\in(0,\infty)$ and condition (\ref{107}) is satisfied for $l=1$ and $g(\cdot)=g_0(\cdot)$.

$(ii)$ Among the above three propositions, it is clear that condition (\ref{103})
is the weakest and condition (\ref{107}) is the strongest.
In Proposition \ref{pron106}, if $EZ_1^r<\infty$ for some $r\ge1$,
then $g(x)=x^{r-1}g_0(x)$ for $x\ge0$,
where $g_0(\cdot)$ is a positive function such that $g_0(x)\underline{\uparrow}\infty$
and $g_0(x)x^{-\omega}\overline{\downarrow}0$ for any $\omega\in(0,1)$.
Therefore, the larger the value of $r$ is, the larger the range of $g(\cdot)$ is;
when $r=1$, the range of $g(\cdot)$ is small.
\end{remark}

\subsection{\normalsize\bf Main results}

In this subsection, we present the key findings of this paper.
These findings encompass the precise large deviations of the finite-time ruin probability and the random sums of net loss process,
as well as the asymptotics of the random-time ruin probability.
These results are based on two dependent structures in the nonstandard renewal risk model mentioned above.

Case 1. $\{Y_i,i\geq1\}$, $\{Z_i,i\geq1\}$ and $\tau$ are independent of each other.
but at least one of $\{Y_i,i\geq1\}$ and $\{Z_i,i\geq1\}$ is a dependent random variable sequence.

%Case 2. $\{(Y_i,Z_i),i\geq1\}$ and $\tau$ are independent of each other,
%$\{(Y_i,Z_i),i\geq1\}$ are independent and identically distributed random variables,
%but $Y_i$ and $Z_i$ are allowed to be arbitrarily dependent for all $i\ge1$.

Case 2. $\{(Y_i,Z_i),i\geq1\}$ and $\tau$ are independent of each other,
but $\{Y_i,i\geq1\}$ and $\{Z_i,i\geq1\}$ are allowed to be arbitrarily dependent.

Clearly, Case 1 is properly contained by Case 2.

We then denote the distributions of $Y_1+cZ_1$ and $\max\{Y_1,cZ_i\}$ as $G_1$ and $G_2$
, respectively, with corresponding means $\mu_{G_1}$ and $\mu_{G_2}$. These distributions and their means will be utilized in the subsequent results and their proofs.

\begin{thm}\label{thm503}
In Case 1 of the above nonstandard renewal risk model, let $Y_i,i\ge1$ be ENOD random variables
with common distribution $G\in\mathcal{C}$ and dominating coefficient $M_G$,
and let $Z_i,\ i\ge1$ be independent and identically distributed random variables with common distribution $H$ satisfying
\begin{eqnarray}\label{532}
\overline{H}(x)=o\big(\overline{G}(x)\big).
\end{eqnarray}

$(i)$ For each $\gamma>\max\{1,\mu_G\}$,%\mu_H^{-1}\mu_{G}
\begin{eqnarray}\label{509000}
a\le\lim\limits_{t\to\infty}\inf\limits_{x\ge\gamma t}\frac{\psi(x;t)}{\mu_H^{-1}t\overline{G}(x-\mu_{G_2}\mu_H^{-1}t+ct)}
\le\lim\limits_{t\to\infty}\sup\limits_{x\ge\gamma t}\frac{\psi(x;t)}{\mu_H^{-1}t\overline{G}(x-\mu_{G_2}\mu_H^{-1}t+ct)}\le1,
\end{eqnarray}
where $a=\overline{G_*}\big(1+(\mu_{G_2}-\mu_G)\mu_H^{-1}(\gamma-\mu_H^{-1}\mu_{G_2}+c)^{-1}\big)$, and
\begin{eqnarray}\label{5090000}
1\le\lim\limits_{t\to\infty}\inf\limits_{x\ge\gamma t}\frac{\psi(x;t)}{\mu_H^{-1}t\overline{G}(x-\mu_{G}\mu_H^{-1}t+ct)}
\le\lim\limits_{t\to\infty}\sup\limits_{x\ge\gamma t}\frac{\psi(x;t)}
{\mu_H^{-1}t\overline{G}(x-\mu_{G}\mu_H^{-1}t+ct)}\le b^{-1},
\end{eqnarray}
where $b=\max\big\{\overline{G_*}\big(1+\frac{\mu_H^{-1}(\mu_{G_2}-\mu_G)}{\gamma-\mu_H^{-1}\mu_{G}+c}\big),\
\overline{G_*}\big(1+\frac{1}{\gamma-1+c-\mu_H^{-1}\mu_{G}}\big)\}$.

$(ii)$  For each $\gamma>0$ and any positive function $s(\cdot)$ satisfying $s(t)\to\infty$ as $t\to\infty$,
\begin{eqnarray}\label{509001}
\lim\limits_{t\to\infty}\sup\limits_{x\ge\gamma ts(t)}\big|\psi(x;t)
\big(\mu_H^{-1}t\overline{G}(x-\mu_{G_2}\ (or\ \mu_G)\ \mu_H^{-1}t+ct)\big)^{-1}-1\big|=0,
\end{eqnarray}
\end{thm}

We can also get the result corresponding to Theorem \ref{thm503} $(i)$ in Case 2,
if we narrow the scope of $x$ in another way.

\begin{thm}\label{thm1020}
In Case 2 of the above nonstandard renewal risk model, let $Y_i,i\ge1$ be ENOD random variables with common distribution $G\in\mathcal{C}$ and dominating coefficient $M_G$,
and let $Z_i,\ i\ge1$ be ENOD random variables with common distribution $H$ satisfying condition (\ref{532})
and dominating coefficient $M_H$.
Then for each pair $\max\{\mu_G,\mu_H^{-1}\mu_{G}\}<\gamma<\Gamma<\infty$, %, $a$ in (\ref{509000}) and $b$ in (\ref{5090000}),
\begin{eqnarray}\label{5090002}
1\le\lim\limits_{t\to\infty}\inf\limits_{x\in[\gamma t,\Gamma t]}\frac{\psi(x;t)}{\mu_H^{-1}t\overline{G}(x-\mu_{G}\mu_H^{-1}t+ct)}
\le\lim\limits_{t\to\infty}\sup\limits_{x\in[\gamma t,\Gamma t]}
\frac{\psi(x;t)}{\mu_H^{-1}t\overline{G}(x-\mu_{G}\mu_H^{-1}t+ct)}\le d^{-1}.
\end{eqnarray}
where $d=\overline{G_*}\big(1+c(\gamma-\mu_H^{-1}\mu_{G})^{-1}\big)$.
\end{thm}

\begin{remark}\label{rem102}
$(i)$ In $\overline{G}(x-\mu_G\mu_H^{-1}t+ct)$ of (\ref{5090000}), (\ref{509001}) and (\ref{5090002}),
$\mu_G\mu_H^{-1}t$ and $ct$ respectively reflect the effects of decentralization
and income $cS_{N(t)}^H$ for $\psi(x;t)$.

$(ii)$ There is no substantive difference between $\lambda(t)$ and $\mu_H^{-1}t$,
because when $t\to\infty$, $\lambda(t)\sim\mu_H^{-1}t$ according to Proposition \ref{pron106}
(thus Proposition \ref{pron105} and Proposition \ref{pron104}).
However, we prefer to replace $\lambda(t)$ with $\mu_H^{-1}t$, which seems more intuitive and explicit.

$(iii)$ From (\ref{33070001}) below, we know that
$$\overline{G_*}^{-1}\Big(1+\frac{\mu_H^{-1}(\mu_{G_2}-\mu_G)}{\gamma-\mu_H^{-1}\mu_{G}+c}\Big)
\le d^{-1}=\overline{G_*}^{-1}\Big(1+\frac{c}{\gamma-\mu_H^{-1}\mu_{G}}\Big).$$
Therefore, we need to find an upper bound $b^{-1}$ smaller than $d^{-1}$ in (\ref{5090000}) by a new method.

$(iv)$ In further research, we hope to get more perfect results,
such as achieving $a=1$ in (\ref{509000}), $b=1$ in (\ref{5090000}) and $c=1$ in (\ref{5090002});
or eliminating the influence of $s(t)$ in (\ref{509001}).
\end{remark}

For net loss process at time $t$ $\big(R_{0}(t)=S^F_{N(t)}:t\ge0\big)$,
we can derive the following two more ideal results than the above two theorems.
Among these results, the one presented in conjunction with Theorem \ref{thm1020} is particularly intriguing.

\begin{thm}\label{thm102}
In the above nonstandard renewal risk model, let $Y_i,i\ge1$ be ENOD random variables
with common distribution $G\in\mathcal{C}$ and dominating coefficients $M_G$,
and let $Z_i,\ i\ge1$ be ENOD random variables with dominating coefficients $M_H$.
%Assume that there exists some almost increasing function $g(\cdot)$ satisfying (\ref{107}) for
%some $l\ge1$ and some $0<\omega<1$ such that $EZ_1g(Z_1)<\infty$.
If (\ref{532}) is satisfied, then the following two conclusions hold.

$(i)$ In Case 1, for each $\gamma>\mu_G$,
\begin{eqnarray}\label{10015}
\lim\limits_{t\to\infty}\sup\limits_{x\ge\gamma t}\Big|P\big(R_{0}(t)>x\big)
\big(\mu_H^{-1}t\overline{G}(x-\mu_G\mu_H^{-1}t+ct)\big)^{-1}-1\Big|=0.
\end{eqnarray}

$(ii)$ In Case 2, for each pair $\mu_G<\gamma<\Gamma<\infty$,
\begin{eqnarray}\label{203}
\lim\limits_{t\to\infty}\sup\limits_{x\in[\gamma t,\Gamma t]}
\Big|P\big(R_{0}(t)>x\big)\big(\mu_H^{-1}t\overline{G}(x-\mu_G\mu_H^{-1}t+ct)\big)^{-1}-1\Big|=0.
\end{eqnarray}
\end{thm}

\begin{remark}\label{rem302}
$(i)$ In Theorem \ref{thm102} $(i)$, if $X_i=Y_i-cZ_i,\ i\ge1$ are ENOD random variables with common distribution $F$,
then the conclusion (\ref{10015}) can also be obtained according to Theorem 1.1 of Tang \cite{T2006} for NOD random variables
and Theorem 2.1 of Liu \cite{L2009} for ENOD random variables under the additional conditions that
$$E(X_1^-)^r<\infty\ \ \ \text{for some}\ r>1\ \ \ \text{and}\ \ \ \ F(-x)=o\big(\overline{F}(x)\big),$$
where $X_1^-=-X_1\emph{\textbf{1}}_{\{X_1<0\}}$.
However, Theorem \ref{thm102} can avoid these additional conditions.

$(ii)$ In Theorem \ref{thm102}, we may can expand the ENOD dependent structure of $Z_i,i\ge1$
to WOD dependent structure under stronger moment conditions.

$(iii)$ In some relevant results, for example, Theorem 1.1 of Chen et al. \cite{CWY2021}, %and Corollary 2.1 of Gao et al. \cite{GLL2023}
requires $\gamma>0$, but Theorem \ref{thm102} $(ii)$ of this paper requires $\gamma>\mu_G$.
In fact, the two conditions are equivalent because the former aims for centralized case and the latter aims for noncentralized case,  see the proof of Lemma \ref{lem201} below.
In addition, the research object of Theorem 1.1 of Chen et al. \cite{CWY2021} is all $S^G_{N(t)}$.
%\textcolor[rgb]{0.00,0.07,1.00}{and that result require that $\{(Y_i,Z_i),i\geq1\}$ are independent and identically distributed random variables}.
\end{remark}

From the proof of Theorem \ref{thm102} $(ii)$, we can immediately find a corresponding result about loss process at time $t$ $\big(R_{00}(t)=\sum_{i=1}^{N(t)}Y_i:t\ge0\big)$,
which improve Theorem 1.1 of Chen et al. \cite{CWY2021}, % and Corollary 2.1 of Gao et al. \cite{GLL2023}
in which $Y_i$ and $Z_i$ are arbitrarily dependent for all $i\ge1$, however $(Y_i,Z_i),i\ge1$ are independent.

\begin{Corol}\label{Corol101}
In Case 2 of the above nonstandard renewal risk model, let $Y_i,i\ge1$ be ENOD random variables
with common distribution $G\in\mathcal{C}$ and dominating coefficients $M_G$.
%Assume that there exists some almost increasing function $g(\cdot)$ satisfying (\ref{107}) for
%some $l\ge1$ and some $0<\omega<1$ such that $EZ_1g(Z_1)<\infty$.
Then for each pair $\mu_G<\gamma<\Gamma<\infty$, %the following two conclusions are hold.
%$(i)$ In Case 1, for each $\gamma>\mu_G$,
%\begin{eqnarray}\label{10015}
%\lim\limits_{t\to\infty}\sup\limits_{x\ge\gamma t}\Big|P\big(R_{00}(t)>x\big)
%\big(\mu_H^{-1}t\overline{G}\big(x-\mu_G\mu_H^{-1}t\big)^{-1}-1\Big|=0.
%\end{eqnarray}
%
%$(ii)$ In Case 2,
\begin{eqnarray}\label{1116}
\lim\limits_{t\to\infty}\sup\limits_{x\in[\gamma t,\Gamma t]}
\Big|P\big(R_{00}(t)>x\big)\big(\mu_H^{-1}t\overline{G}(x-\mu_G\mu_H^{-1}t)\big)^{-1}-1\Big|=0.
\end{eqnarray}
\end{Corol}

%Recall the concept of random-time ruin probability in (\ref{404}).
%For some relevant researches, see Remark \ref{rem101} below.
Finally, we give a new result on asymptotics of the random-time ruin probability,
which is closely related to the above results on the precise large deviations.
%in Case 3,
%where $\{Y_i,i\geq1\}$ and $\{Z_i,i\geq1\}$ are allowed to be arbitrarily dependent.
%Therefore, the result is one of the main results of this paper.

\begin{thm}\label{thm602}
In Case 1 of the above nonstandard renewal risk model,
let $Y_i,i\ge1$ be ENOD random variables with a common distribution $G\in\mathcal{C}$,
and let $Z_i,i\ge1$ be WOD random variables with a common distribution $H$, dominating coefficients $g_{L,H}(n),g_{U,H}(n),n\ge1$.
Assume that there exists some almost increasing function $g(\cdot)$,
some integer $m\ge1$ and some $\omega\in(0,1)$ such that $EZ_1g(Z_1)<\infty$ and condition (\ref{107}) is satisfied.
Further, if (\ref{532}) is satisfied and
\begin{eqnarray}\label{408}
P(\tau>x)=o\big(\overline{G}(x)\big),
\end{eqnarray}
then $\max\{E\tau,EN(\tau)\}<\infty$ and
\begin{eqnarray}\label{409}
\psi(x;\tau)\sim EN(\tau)\overline{G}(x).
\end{eqnarray}
%in Case 1 or in Case 2 and $Z_i,i\ge1$ are ENOD random variables.
\end{thm}

\begin{remark}\label{rem104}
$(i)$ We make a comparison between the above theorem and the existing related results.

In Theorem 2.1 of Wang et al. \cite{WGWL2009} and Corollary 2.2 (2) of Wang et al. \cite{WCWM2012},
the random variables $Y_i,i\ge1$ are independent and identically distributed.
In addition, the former requires $Ee^{s_0N(\tau)}<\infty$ for some $s_0>0$,
and the latter requires $EZ_1^p<\infty$ for some $p\ge2$ or some other specified condition.

In Theorems 3.1, 3.2, 3.3 and 3.4 of Wang et al. \cite{WGWL2009},
the random variables $Y_i,i\ge1$ are NUOD and the random variables $Z_i,i\ge1$ are NOD.
Furthermore, Theorems 3.1 and 3.2 above require $EY_1^{1+\beta}<\infty$ for some $\beta>0$.
Additionally, Theorems 3.3 and 3.4 above require $EN^p(\tau)<\infty$ for some $p>J_G^+$,
where $J_G^+\ge1$. %which is stronger than condition (\ref{408}), see Proposition \ref{pron107} below.

All the above results require that $Z_i,i\ge1$ are independent of $Y_i,i\ge1$.

%In addition, these known results require that $(Y_i:i\ge1)$ and $(Z_i:i\ge1)$ are independent of each other.

$(ii)$ Here, we give a sufficient condition for (\ref{408}) that,
if $G\in\mathcal{D}$ and $E\tau^p<\infty$ for some $p>J_G^+$, then (\ref{408}) holds.
%Clearly, Theorem \ref{thm102} of this paper has more advantages than the existing results mentioned above.
%This is because we mainly make use method of precise large deviation to prove the theorem,
%which is different from that of the above existing results.
%In other words, this result also reflects the significance and function of the study of precise large deviation.
\end{remark}

%In the above theorem, especially, when $\tau=t\ a.s.$ for some $t>0$,
%we immediately obtain the asymptotics of finite-time ruin probability.
%
%\begin{Corol}\label{cor401}
%In Case 3 of the above nonstandard renewal risk model,
%let $Y_i,i\ge1$ be ENOD random variables with a common distribution $G\in\mathcal{C}$,
%and let $Z_i,i\ge1$ be WOD random variables with a common distribution $H$, dominating coefficients $g_{L,H}(n),g_{U,H}(n),n\ge1$.
%Assume that there exists some almost increasing function $g(\cdot)$,
%some integer $m\ge1$ and some $\omega\in(0,1)$ such that $EZ_1g(Z_1)<\infty$ and condition (\ref{107}) is satisfied.
%Further, if (\ref{532}) is satisfied, then for each $t>0$,
%\begin{eqnarray}\label{437}
%\psi(x;t)\sim \mu_H^{-1}t\overline{G}(x).
%\end{eqnarray}
%\end{Corol}

The paper is organized as follows.
We prove these results in Section 3.
To this end, we introduce a preliminary result on the precise large deviations for sums
of nonnegative and noncentralized ENOD random variables in Section 2.
Then, in section 4, we give some applications of the above results including precise large deviations
for the proportional-net-loss process and the excess-of-net-loss process,
and the asymptotic estimates of the mean of stop-net-loss reinsurance treaty.

\section{\normalsize\bf A lemma}
\setcounter{equation}{0}\setcounter{thm}{0}\setcounter{lemma}{0}\setcounter{remark}{0}\setcounter{exam}{0}

The earlier work on precise large deviations for nonrandom sums of independent
and identically distributed random variables with a common distribution $F$
can be found in Heyde $\cite{H1967a},\cite{H1967b},\cite{H1968}$ and Nagaev $\cite{N1969a},\cite{N1969b}$.
With the development of research on precise large deviations,
the classic case of the distribution $F\in\mathcal{R}_{\alpha}$ on $[0,\infty)$ is attributed to Nagaev $\cite{N1973},\cite{N1979}$;
the case of $F\in\mathcal{ERV}(\alpha,\beta)$ can be referred to Cline and Hsing $\cite{CH1991}$,
Kl\"{u}ppelberg and Mikosch $\cite{KM1997}$, Mikosch and Nagaev $\cite{MN1998}$, Tang et al. $\cite{TSJZ2001}$, among others;
and the case of $F\in\mathcal{C}$ is shown in Ng et al. $\cite{NTYY2004}$ and Wang and Wang $\cite{WW2007}$.
For the corresponding results of local probability and density, see Doney \cite{D1989}, Denisov et al. \cite{DDS2008}, Lin \cite{L2008},
Yang et al. \cite{YLS2010}, Cheng and Li \cite{CL2016}, etc.

For dependent random variables, there are numerous corresponding findings.
Konstantinides and Mikosch $\cite{KM2005}$ obtained several results
regarding the precise large deviations of sums for a stochastic recurrence equation,
Mikosch and Wintenberger $\cite{MW2013}$ delved into this topic for a stationary regularly varying sequence of random variables,
while, Mikosch and Rodionov \cite{MR2021} further explored the subject for the partial sums of a stationary sequence with a subexponential
marginal distribution, such as regularly varying or lognormal distributions.
Among them, the results about lognormal distribution are particularly striking.

In this paper, we focus on another kind of dependent structure of random variables.
Theorem 3.1 of Tang $\cite{T2006}$ studied the precise large deviation of NOD random variables
with common distribution $F\in\mathcal{C}$ on $(-\infty,\infty)$.
Thereafter, Theorem 2.1 of Liu $\cite{L2009}$ replaced the condition $xF(-x)=o\big(\overline{F}(x)\big)$ in Tang $\cite{T2006}$ with the condition $F(-x)=o\big(\overline{F}(x)\big)$ and gave the same result for ENOD random variables.
Kl\"{u}ppelberg and Mikosch $\cite{KM1997}$, Chen et al. \cite{CYN2011} and Gao et al. \cite{GLL2023} gave some corresponding results for random sums of independent and ENOD random variables, respectively.
In addition, Theorem 2.1 of Wang et al. $\cite{WWC2006}$ studied the precise large deviation for $F\in\mathcal{D}$ on $[0,\infty)$.
%Generally, Theorem 1 of Wang et al. \cite{WYL2012} discusses the case of WOD random variables with different distributions.
%However, in Wang et al. \cite{WYL2012}, the range of the dominating coefficients of WOD random variables is too small
%to obtain the corresponding result as the following (\ref{106000}) and (\ref{Theorem1030}).
%In our view, Corollary 1 of Wang et al. \cite{WYL2012} gives such a result for identically distributed WOD random variables that, for each integer $m\ge1$ and $v\in\big(0,m(m+1)^{-1}\big)$, if $\max\{g_U(n),\ g_L(n)\}=O(n^{m^{-1}})$ for all $n\ge1$, then for each $\gamma>0$, $$\overline{F}_*(v^{-1})\le\liminf\limits_{n\to\infty}\inf\limits_{x\geq\gamma n}\frac{P(S_n>x)}{n\overline{F}(x)} \le\limsup\limits_{n\to\infty}\sup\limits_{x\geq\gamma n}\frac{P(S_n>x)}{n\overline{F}(x)}\leq\overline{F}_*^{-1}(v^{-1}).$$ Clearly, when $v\uparrow1$, $\overline{F}_*(v^{-1})\uparrow L_F$ and $m\uparrow\infty$, the latter implies that $X_i,\ i\ge1$ be ENOD.

Most of the aforementioned results necessitate a mean value of zero for random variables.
However, in order to establish the majority of the theorems in this paper,
it is essential to convert some results into a noncentralized form.

The following result is a noncentralized version of Theorem 2.2 of Liu $\cite{L2009}$
with a direct proof for nonnegative random variables.

\begin{lemma}\label{lem201}
In the above nonstandard renewal risk model, let $Y_i,i\geq 1$ be ENOD random variables
with common distribution $G\in\mathcal{C}$ on $[0,\infty)$, dominating coefficient $M_G$ and mean $\mu_G$.
Then for each $\gamma>\mu_G$,
\begin{eqnarray}\label{10600}
\lim\limits_{n\to\infty}\sup\limits_{x\in[\gamma n,\infty)}\big|P(S^G_n>x)\big(n\overline{G}(x-\mu_G n)\big)^{-1}-1\big|=0.
\end{eqnarray}
\end{lemma}
\proof Let $X_1=Y_1-\mu_G$ with distribution $F$.
Clearly, $EX_1=0$, $P(X_1^->x)=o\big(\overline{F}(x)\big)$, $E(X_1^-)^r<\infty$ for each $r>1$ and
$$P(S^G_n>x)=P\Big(\sum_{i=1}^n(Y_i-\mu_G)>x-n\mu_G\Big)=P\Big(\sum_{i=1}^nX_i>x-n\mu_G\Big).$$
Further, for each $\gamma>\mu_G$, by $x\ge\gamma n$, we have $\gamma_0=\gamma-\mu_G>0$ and
$$x-n\mu_G\ge(\gamma-\mu_G)n=\gamma_0n.$$
In fact, $\gamma_0>0$ is equivalent to $\gamma>\mu_G$.
Therefore, according to Theorem 2.2 of Liu \cite{L2009},
\begin{eqnarray*}
\lim\limits_{n\to\infty}\sup\limits_{x\ge\gamma n}\Big|\frac{P(S^G_n>x)}{n\overline{G}(x-\mu_G n)}-1\Big|
=\lim\limits_{n\to\infty}\sup\limits_{x-n\mu_G\ge\gamma_0 n}\Big|\frac{P(S^G_n-n\mu_G>x-n\mu_G)}{n\overline{G}(x-\mu_G n)}-1\Big|=0,
\end{eqnarray*}
that is (\ref{lem201}) holds.
$\hspace{\fill}\Box$
\begin{remark}\label{rem201}
%The lemma plays a key role in proofs of Theorem \ref{thm101} and Theorem \ref{thm103}. And it is not difficult to find that
In the lemma, the restrict $\gamma>\mu_G$ has no effect on some applications, %of the in some cases,
see, for example, %on asymptotic behavior of random time ruin probability in
the proof of Theorem \ref{thm602}.
\end{remark}

\section{\normalsize\bf Proofs of main results}
\setcounter{equation}{0}\setcounter{thm}{0}\setcounter{pron}{0}\setcounter{lemma}{0}\setcounter{Corol}{0}

\subsection{\normalsize\bf Proof of Theorem \ref{thm503}}

$(i)$ For some integer $k_0\ge\max\{2, 2\mu_H^{-1}\}$, we take any integer $k\ge k_0$ and any $\delta\in(0,1)$ small enough to split
\begin{eqnarray}\label{3303}
&&\psi(x;t)=\bigg(\sum_{m=1}^{\lfloor\frac{(1+\delta)t}{k\mu_H}\rfloor}
+\sum_{m=\lfloor\frac{(1+\delta)t}{k\mu_H}\rfloor+1}^{\lfloor\frac{(1+\delta)t}{\mu_H}\rfloor}
+\sum_{m=\lfloor\frac{(1+\delta)t}{\mu_H}\rfloor+1}^{\infty}\bigg)
P\Big(\max_{1\le n\le m}\sum_{i=1}^{n}(Y_i-cZ_i)>x,N(t)=m\Big)\nonumber\\
&=&\psi_1(x;t)+\psi_2(x;t)+\psi_3(x;t).
\end{eqnarray}

For $\psi_1(x;t)$, by $x\ge\gamma t$ for each $\gamma>\mu_G$, we have
$$x\ge\gamma t>\mu_Gt=\mu_G(1+\delta)^{-1}k\mu_H(1+\delta)(k\mu_H)^{-1}t=\gamma_1(1+\delta)(k\mu_H)^{-1}t\ \ \ \text{
and}\ \ \ \gamma_1>\mu_G.$$
Therefore, according to Lemma \ref{lem201}, by $\mu_G<c\mu_H$, when $t\to\infty$, it holds uniformly for all $x\ge\gamma t$ that
\begin{eqnarray}\label{3304}
&&\psi_1(x;t)\le P\bigg(\sum_{i=1}^{\lfloor(1+\delta)(k\mu_H)^{-1}t\rfloor}Y_i>x\bigg)\nonumber\\
&\lesssim&(1+\delta)(k\mu_H)^{-1}t\overline{G}\big(x-(1+\delta)k^{-1}\mu_G\mu_H^{-1}t\big)\nonumber\\
&<&2k^{-1}\mu_H^{-1}t\overline{G}\big(x-\mu_G\mu_H^{-1}t+ct-(1+\delta)k^{-1}\mu_G\mu_H^{-1}t+\mu_G\mu_H^{-1}t-ct\big)\nonumber\\
&\le&2k^{-1}\mu_H^{-1}t\overline{G}\Big((x-\mu_G\mu_H^{-1}t+ct)\big(1-(ct-(1-2k^{-1})\mu_G\mu_H^{-1}t)(x-\mu_G\mu_H^{-1}t+ct)^{-1}\big)\Big)\nonumber\\
&\le&2k^{-1}\mu_H^{-1}t\overline{G}\Big((x-\mu_G\mu_H^{-1}t+ct)\big(1-(c-(1-2k^{-1})\mu_G\mu_H^{-1})(\gamma-\mu_G\mu_H^{-1}+c)^{-1}\big)\Big)\nonumber\\
&=&k^{-1}O\big(\mu_H^{-1}t\overline{G}(x-\mu_G\mu_H^{-1}t+ct)\big)\nonumber\\
&=&o\big(\mu_H^{-1}t\overline{G}(x-\mu_G\mu_H^{-1}t+ct)\big)\ \ \ \ \ \ \text{then let}\ \ k\to\infty.
\end{eqnarray}

To deal with $\psi_2(x;t)$, we first give two lemmas.

\begin{lemma}\label{lem3302}
Let $\xi_i,i\ge1$ be WUOD (or WLOD) random variables with common distribution $V_1$ on $(-\infty,\infty)$
and dominating coefficients $g_{U,V_1}(n)$ \big(or $g_{L,V_1}(n)\big),\ n\ge1$,
and let $\eta_i,i\ge1$ be independent random variables with common distribution $V_2$ on $(-\infty,\infty)$.
If $\xi_i,i\ge1$ and $\eta_i,i\ge1$ are independent of each other,
then $\xi_i+\eta_i,i\ge1$ still is WUOD (or WLOD) random variables
with common distribution $W$ on $(-\infty,\infty)$ and the same dominant coefficients as $\xi_i,i\ge1$.
\end{lemma}

\begin{lemma}\label{lem3301}
Let $\xi_1$ and $\eta_1$ be two random variables that are allowed to be arbitrarily dependent
with their respective distributions $V_1$ and $V_2$ on $[0,\infty)$.
Let $W_1$ be the distribution of $\xi_1+\eta_1$ and $W_2$ be that of $\max\{\xi_1,\eta_1\}$.

$(i)$ If $V_1\in\mathcal{C}$ and (\ref{532}) is satisfied with $G=V_1$ and $H=V_2$, then $W_i\in\mathcal{C}$ and
\begin{eqnarray}\label{5336}
\overline{W_i}(x)\sim\overline{V_1}(x)\ \ \ \ \text{for}\ \ i=1,2.
\end{eqnarray}

$(ii)$ Further, let $\xi_i,i\ge1$ be WUOD (or WLOD) random variables with common distribution $V_1$ and dominating coefficients
$g_{U,V_1}(n)$\ \big(or $g_{L,V_1}(n)\big)$ for all $n\ge1$,
and let $\eta_i,i\ge1$ be independent random variable with common distribution $V_2$.
If $\xi_i,i\ge1$ and $\eta_i,i\ge1$ are independent of each other,
then $\max\{\xi_i,\eta_i\},i\ge1$ also are WUOD (or WLOD) random variables.
\end{lemma}

\proof $(i)$ If $\xi$ is independent of $\eta$, then the conclusion is known,
see, for example, Theorem 1.1 of Cheng et al \cite{CNPW2012}.
Generally, by $V_1\in\mathcal{C}$, (\ref{532}) and
$$\overline{V_1}(x)\le\overline{W_2}(x)\le\overline{W_1}(x)\le\overline{V_1}\big((1-\delta)x\big)+\overline{V_2}(\delta x)$$
for any $0<\delta<1$, we immediately reach the above conclusion.

$(ii)$ If the random variables $\xi_i,i\ge1$ are WUOD, then according to Lemma \ref{lem3302},  the random variables $\xi_i+\eta_i,i\ge1$ are also WUOD. This, in turn, implies that the random variables $\max\{\xi_i,\eta_i\},i\ge1$ are WUOD as well.
In fact, for each $n\ge1$ and $x_i\in(-\infty,\infty),\ 1\le i\le n$, we have
\begin{eqnarray}\label{2101}
&&P\Big(\bigcap^{n}_{i=1}\{\max\{\xi_i,\eta_i\}>x_i\}\Big)\leq P\Big(\bigcap^{n}_{i=1}\{\xi_i+\eta_i>x_i\}\Big)\nonumber\\
&\le&g_{U,V_1}(n)\prod_{i=1}^n P(\xi_i+\eta_i> x_i)\nonumber\\
&\le&g_{U,V_1}(n)\prod_{i=1}^n P(2\max\{\xi_i,\eta_i\}> x_i)\nonumber\\
&=&O\Big(\prod_{i=1}^n P(\max\{\xi_i,\eta_i\}> x_i)\Big).
\end{eqnarray}
If $\xi_i,i\ge1$ are WLOD random variables, then by
\begin{eqnarray}\label{2102}
&&P\Big(\bigcap^{n}_{i=1}\{\max\{\xi_i,\eta_i\}\le x_i\}\Big)=P\Big(\bigcap^{n}_{i=1}\{\xi_i\le x_i,\eta_i\le x_i\}\Big)\nonumber\\
&\le&g_{L,V_1}(n)\prod_{i=1}^n P(\xi_i\le x_i)P(\eta_i\le x_i)\nonumber\\
&=&g_{L,V_1}(n)\prod_{i=1}^n P(\max\{\xi_i,\eta_i\}\le x_i),
\end{eqnarray}
we know that $\max\{\xi_i,\eta_i\},i\ge1$ still are WLOD random variables.
$\hspace{\fill}\Box$\\

In Case 1, according to Lemma \ref{lem3302} and Lemma \ref{lem3301}, $Y_i+cZ_i,i\ge1$ and $\max\{Y_i,cZ_i\},i\ge1$ still are ENOD random variables with common distribution $G_1$ and $G_2$, respectively.
In addition, $G_i\in\mathcal{C}$ and $\overline{G_i}(x)\sim\overline{G}(x),\ i=1,2$.

Now we deal with $\psi_2(x;t)$.

For each $\gamma>1$, by $c>\mu_G\mu_H^{-1}$, we have
$$\big((\gamma+c)\mu_H-\mu_G-\mu_H\big)2^{-1}(\mu_G+\mu_H)^{-1}>0.$$
Then we take $0<\delta<\big((\gamma+c)\mu_H-\mu_G-\mu_H\big)2^{-1}(\mu_G+\mu_H)^{-1}$ such that, when $t\ge\delta^{-1}\mu_H$,
\begin{eqnarray*}
&&x+ct\ge(\gamma+c)t=(\gamma+c)t\big((1+\delta)\mu_H^{-1}t+1\big)^{-1}\big((1+\delta)\mu_H^{-1}t+1\big)\nonumber\\
&\ge&(\gamma+c)\big((1+2\delta)\mu_H^{-1}\big)^{-1}\big((1+\delta)\mu_H^{-1}t+1\big)=\gamma_2\big((1+\delta)\mu_H^{-1}t+1\big)
\end{eqnarray*}
and $\gamma_2>\mu_G+\mu_H>\mu_{G_2}$.
Therefore, according to Lemma \ref{lem201}, by $\mu_G<c\mu_H$ and $\overline{G_2}(x)\sim\overline{G}(x)$,
for the above $\delta$ and each $\gamma>1$, when $t\to\infty$, it holds uniformly for all $x\ge \gamma t$ that
\begin{eqnarray}\label{33070}
&&\psi_2(x;t)=P\bigg(\bigcup_{m=\lfloor(1+\delta)(k\mu_H)^{-1}t\rfloor+1}^{\lfloor (1+\delta)\mu_H^{-1}t\rfloor}
\Big\{\max_{1\le n\le m}\Big(S^G_n-cS^H_{m+1}+\sum_{i=n+1}^{m+1}cZ_i\Big)>x,N(t)=m\Big\}\bigg)\nonumber\\
&\le&P\bigg(\bigcup_{m=\lfloor(1+\delta)(k\mu_H)^{-1}t\rfloor+1}^{\lfloor (1+\delta)\mu_H^{-1}t\rfloor}
\Big\{\max_{1\le n\le m}\Big(S^G_n+\sum_{i=n+1}^{m+1}cZ_i\Big)>x+ct,N(t)=m\Big\}\bigg)\nonumber\\\
&\le&P\bigg(\max_{1\le n\le\lfloor(1+\delta)\mu_H^{-1}t\rfloor}
\Big(S^G_n+\sum_{i=n+1}^{\lfloor(1+\delta)\mu_H^{-1}t\rfloor+1}cZ_i\Big)>x+ct,
\bigcup_{m=\lfloor(1+\delta)(k\mu_H)^{-1}t\rfloor+1}^{\lfloor (1+\delta)\mu_H^{-1}t\rfloor}\{N(t)=m\}\bigg)\nonumber\\
&\le&P\Big(S^{G_2}_{\lfloor(1+\delta)\mu_H^{-1}t\rfloor+1}>x+ct\Big)\nonumber\\
&\sim&(1+\delta)\mu_H^{-1}t\overline{G_2}(x-(1+\delta)\mu_H^{-1}\mu_{G_2}t+ct)\nonumber\\
&=&(1+\delta)\mu_H^{-1}t\overline{G_2}\Big((x-\mu_H^{-1}\mu_{G_2}t+ct)
\big(1-\delta\mu_{G_2}\mu_H^{-1}t(x-\mu_H^{-1}\mu_{G_2}t+ct)^{-1}\big)\Big)\nonumber\\
&\le&(1+\delta)\mu_H^{-1}t\overline{G_2}\Big((x-\mu_H^{-1}\mu_{G_2}t+ct)
\big(1-\delta\mu_{G_2}\mu_H^{-1}(\gamma-\mu_H^{-1}\mu_{G_2}+c)^{-1}\big)\Big)\nonumber\\
&\lesssim&(1+\delta)\mu_H^{-1}t\overline{G}(x-\mu_H^{-1}\mu_{G_2}t+ct)
\overline{G_*}^{-1}\Big(\big(1-\delta\mu_{G_2}\mu_H^{-1}(\gamma-\mu_H^{-1}\mu_{G_2}+c)^{-1}\big)^{-1}\Big)\nonumber\\
&\sim&\mu_H^{-1}t\overline{G}(x-\mu_H^{-1}\mu_{G_2}t+ct)\ \ \ \ \ \ \text{then let}\ \  \delta\downarrow0.
\end{eqnarray}

If we want to replace $\overline{G}(x-\mu_H^{-1}\mu_{G_2}t+ct)$ with $\overline{G}(x-\mu_H^{-1}\mu_{G}t+ct)$,
we have the following conclusions.

From (\ref{33070}), when $t\to\infty$, it holds uniformly for all $x\ge \gamma t$ that
\begin{eqnarray}\label{3307001}
&&\psi_2(x;t)\lesssim\mu_H^{-1}t\overline{G}\big((x-\mu_H^{-1}\mu_{G}t+ct\big)
\big(1-\mu_H^{-1}(\mu_{G_2}-\mu_G)(\gamma-\mu_H^{-1}\mu_{G}+c)^{-1}\big)\Big)\nonumber\\
&\lesssim&\mu_H^{-1}t\overline{G}(x-\mu_H^{-1}\mu_{G}t+ct)
\overline{G_*}^{-1}\big(1+\mu_H^{-1}(\mu_{G_2}-\mu_G)(\gamma-\mu_H^{-1}\mu_{G}+c)^{-1}\big).
\end{eqnarray}

Similarly, for each $\gamma>1$, by $c>\mu_G\mu_H^{-1}$, we take the above $\delta$ and $t$ such that
\begin{eqnarray*}
x+ct\ge\gamma_2\big((1+\delta)\mu_H^{-1}t+1\big)\ \ \ \text{and}\ \ \ \gamma_2>\mu_G+\mu_H=\mu_{G_1}.
\end{eqnarray*}
Therefore, according to Lemma \ref{lem201}, by $\mu_G<c\mu_H$ and $\overline{G_1}(x)\sim\overline{G}(x)$,
for the above $\delta$ and each $\gamma>1$, when $t\to\infty$, it holds uniformly for all $x\ge \gamma t$ that
\begin{eqnarray}\label{330700}
&&\psi_2(x;t)\le P\Big(S^{G_1}_{\lfloor(1+\delta)\mu_H^{-1}t\rfloor+1}>x+ct\Big)\nonumber\\
&\sim&(1+\delta)\mu_H^{-1}t\overline{G_1}\big(x-(1+\delta)\mu_H^{-1}(\mu_{G}+\mu_H)t+ct\big)\nonumber\\
&\lesssim&(1+\delta)\mu_H^{-1}t\overline{G_1}\big((x-\mu_H^{-1}\mu_{G}t+ct\big)
\big(1-(1+\delta+\delta\mu_G\mu_H^{-1})(\gamma-\mu_H^{-1}\mu_{G}+c)^{-1}\big)\Big)\nonumber\\
&\lesssim&(1+\delta)\mu_H^{-1}t\overline{G}(x-\mu_H^{-1}\mu_{G}t+ct)
\overline{G_*}^{-1}\Big(\big(1-(1+\delta+\delta\mu_G\mu_H^{-1})(\gamma-\mu_H^{-1}\mu_{G}+c)^{-1}\big)^{-1}\Big)\nonumber\\
&\sim&\mu_H^{-1}t\overline{G}(x-\mu_H^{-1}\mu_{G}t+ct)
\overline{G_*}^{-1}\big(1+(\gamma-1+c-\mu_H^{-1}\mu_{G})^{-1}\big),\ \ \ \ \text{then let}\ \ \ \delta\downarrow0.
\end{eqnarray}

Or, according to Lemma \ref{lem201}, by $\mu_G<c\mu_H$ and (\ref{3303}), for each $\gamma>\mu_H^{-1}\mu_G,$
when $t\to\infty$, it holds uniformly for all $x\ge\gamma t$ that
\begin{eqnarray}\label{3307000}
&&\psi_2(x;t)\le P(S^{G}_{\lfloor(1+\delta)\mu_H^{-1}t\rfloor+1}>x)\nonumber\\
&\sim&(1+\delta)\mu_H^{-1}t\overline{G}\big(x-(1+\delta)\mu_H^{-1}\mu_{G}t\big)\nonumber\\
&\le&(1+\delta)\mu_H^{-1}t\overline{G}\big((x-\mu_H^{-1}\mu_{G}t+ct\big)
\big(1-(c+\delta\mu_H^{-1}\mu_G)(\gamma-\mu_H^{-1}\mu_{G}+c)^{-1}\big)\Big)\nonumber\\
&\lesssim&(1+\delta)\mu_H^{-1}t\overline{G}(x-\mu_H^{-1}\mu_{G}t+ct)
\overline{G_*}^{-1}\Big(\big(1-(c+\delta\mu_G\mu_H^{-1})(\gamma-\mu_H^{-1}\mu_{G}+c)^{-1}\big)^{-1}\Big)\nonumber\\
&\sim&\mu_H^{-1}t\overline{G}(x-\mu_H^{-1}\mu_{G}t+ct)
\overline{G_*}^{-1}\big(1+c(\gamma-\mu_H^{-1}\mu_{G})^{-1}\big),\ \ \ \ \text{then let}\ \ \ \delta\downarrow0.
\end{eqnarray}

However, for $\overline{G_*}^{-1}\big(1+\frac{\mu_H^{-1}(\mu_{G_2}-\mu_G)}{\gamma-\mu_H^{-1}\mu_{G}+c}\big)$ in (\ref{3307001})
and $\overline{G_*}^{-1}\big(1+\frac{c}{\gamma-\mu_H^{-1}\mu_{G}}\big)$ in (\ref{3307000}), by
$$\mu_G\le\mu_{G_2}\le\mu_{G_1}=\mu_G+c\mu_H,$$
we know that
\begin{eqnarray}\label{33070001}
\overline{G_*}^{-1}\Big(1+\frac{\mu_H^{-1}(\mu_{G_2}-\mu_G)}{\gamma-\mu_H^{-1}\mu_{G}+c}\Big)
\le \overline{G_*}^{-1}\Big(1+\frac{c}{\gamma-\mu_H^{-1}\mu_{G}+c}\Big)
\le\overline{G_*}^{-1}\Big(1+\frac{c}{\gamma-\mu_H^{-1}\mu_{G}}\Big).
\end{eqnarray}
And for $\overline{G_*}^{-1}\big(1+\frac{c}{\gamma-\mu_H^{-1}\mu_{G}}\big)$ in (\ref{3307000})
and $\overline{G_*}^{-1}\big(1+\frac{1}{\gamma-\mu_H^{-1}\mu_{G}+c-1}\big)$ in (\ref{330700}),
when $c\in(0,1]$, we have
\begin{eqnarray}\label{33070002}
\overline{G_*}^{-1}\Big(1+\frac{c}{\gamma-\mu_H^{-1}\mu_{G}}\Big)
\le \overline{G_*}^{-1}\Big(1+\frac{1}{\gamma-\mu_H^{-1}\mu_{G}}\Big)
\le\overline{G_*}^{-1}\Big(1+\frac{1}{\gamma-\mu_H^{-1}\mu_{G}+c-1}\Big).
\end{eqnarray}
When $c>1$, if $\mu_H^{-1}(\mu_{G_2}-\mu_G)\le1$, then
\begin{eqnarray}\label{33070003}
\overline{G_*}^{-1}\Big(1+\frac{\mu_H^{-1}(\mu_{G_2}-\mu_G)}{\gamma-\mu_H^{-1}\mu_{G}+c}\Big)
\le\overline{G_*}^{-1}\Big(1+\frac{1}{\gamma-\mu_H^{-1}\mu_{G}+c-1}\Big);
\end{eqnarray}
if $\mu_H^{-1}(\mu_{G_2}-\mu_G)>1$, then for $\gamma$ large enough,
\begin{eqnarray}\label{33070004}
\overline{G_*}^{-1}\Big(1+\frac{\mu_H^{-1}(\mu_{G_2}-\mu_G)}{\gamma-\mu_H^{-1}\mu_{G}+c}\Big)
\ge\overline{G_*}^{-1}\Big(1+\frac{1}{\gamma-\mu_H^{-1}\mu_{G}+c-1}\Big).
\end{eqnarray}

Therefore, among the three asymptotic upper bounds in (\ref{3307001}), (\ref{330700}) and (\ref{3307000}), we should abandon the third one.

Finally, we deal with $\psi_3(x;t)$. To this end, we first recall Lemma 2.1 of Wang et al. \cite{WWC2006} in the case of NA random variables, which can also easily be generalized to the case of WUOD random variable.

\begin{lemma}\label{lem304}
Let $Y_i,i\ge1$ be WUOD nonnegative random variables with common distribution $G$,
mean $\mu_G\in(0,\infty)$ and dominating coefficients $g_{U,G}(n),\ n\ge1$.
Then, for each $n\ge1$ and any $v>0$,
\begin{eqnarray}\label{4004}
P(S^G_n>x)\le n\overline{G}(v^{-1}x)+g_{U,G}(n)(e\mu_Gnx^{-1})^v,\ \ \ \ \ x\in(0,\infty).
\end{eqnarray}
\end{lemma}
Then according to Proposition \ref{Proposition102}, Proposition \ref{pron106}
(thus Proposition \ref{pron104} and Proposition \ref{pron105}),
Lemma \ref{lem201} and H\"{o}lder inequality,
by $c\mu_H>\mu_G$ and (\ref{4004}) with $g_{U,G}(m)=M_G,\ v>J_G^+,\ J_G^+\ge1$,
when $t\to\infty$, for each $\gamma>\mu_G$, it holds uniformly for all $x\ge\gamma t$ that
\begin{eqnarray}\label{3312}
&&\psi_3(x,t)\le\sum_{m>(1+\delta)\mu^{-1}_Ht}P(S_m^G>x)P\big(N(t)=m\big)\nonumber\\
&\le&\sum_{m>(1+\delta)\mu^{-1}_Ht}\big(m\overline{G}(v^{-1}x)+M_G(e\mu_Gmx^{-1})^v\big)P\big(N(t)=m\big)\nonumber\\
&=&\overline{G}(v^{-1}x)EN(t)\textbf{1}_{\{N(t)>(1+\delta)\mu_H^{-1}t\}}
+M_G(e\mu_Gx^{-1})^vEN(t)N^{v-1}(t)\textbf{1}_{\{N(t)>(1+\delta)\mu_H^{-1}t\}}\nonumber\\
&=&O\Big(\overline{G}(x)E^{2^{-1}}N^2(t)\big(P^{2^{-1}}(N(t)>(1+\delta)\mu_H^{-1}t)
+E^{2^{-1}}N^{2(v-1)}(t)\textbf{1}_{\{N(t)>(1+\delta)\mu_H^{-1}t\}}\big)\Big)\nonumber\\
&=&o\Big(\mu_H^{-1}t\overline{G}\big((x-\mu_G\mu_H^{-1}t+ct)\big(1-(ct-\mu_G\mu_H^{-1}t)(x-\mu_G\mu_H^{-1}t+ct)^{-1}
\big)\big)\Big)\nonumber\\
&=&o\Big(\mu_H^{-1}t\overline{G}\big((x-\mu_G\mu_H^{-1}t+ct)\big(1-(c-\mu_G\mu_H^{-1})(\gamma-\mu_G\mu_H^{-1}+c)^{-1}\big)\big)\Big)\nonumber\\
&=&o\big(\mu_H^{-1}t\overline{G}(x-\mu_G\mu_H^{-1}t+ct)\big).
\end{eqnarray}

Combining (\ref{3303}), (\ref{3304}), (\ref{33070})-(\ref{33070004}) and (\ref{3312}), we know that,
for each $\gamma>\max\{1,\mu_G\}$,
\begin{eqnarray}\label{50901}
\limsup\limits_{t\to\infty}\sup\limits_{x\ge\gamma t}
\psi(x;t)\big(\mu_H^{-1}t\overline{G}(x-\mu_{G_2}\mu_H^{-1}t+ct)\big)^{-1}\le 1.
\end{eqnarray}
and
\begin{eqnarray}\label{50900}
\limsup\limits_{t\to\infty}\sup\limits_{x\ge\gamma t}
\psi(x;t)\big(\mu_H^{-1}t\overline{G}(x-\mu_{G}\mu_H^{-1}t+ct)\big)^{-1}\le b^{-1}.
\end{eqnarray}

On the other hand, for each $\gamma>0$ and any $0<\delta<2^{-1}$, by $c\mu_H>\mu_G$, we have
$$x+ct>(\gamma+c)(1-\delta)^{-1}\mu_H(1-\delta)\mu_H^{-1}t=\gamma_3(1-\delta)\mu_H^{-1}t$$
and
$$\gamma_3>c\mu_H>\mu_G.$$
Thus, according to Lemma \ref{lem201} and Proposition \ref{pron106}, by $\sum_{i=1}^{N(t)}Z_i\le t$,
we know that, when $t\to\infty$, for each $\gamma>0$, it holds uniformly for $x\ge\gamma t$ that
\begin{eqnarray}\label{53351}
&&\psi(x;t)\ge P\Big(\sum_{i=1}^{N(t)}(Y_i-cZ_i)>x\Big)\nonumber\\
&\ge&P\Big(\sum_{i=1}^{N(t)}Y_i>x+ct,(1-\delta)\mu_H^{-1}t\le N(t)\le(1+\delta)\mu_H^{-1}t\Big)\nonumber\\
&\ge&P\big(S^G_{(1-\delta)\mu_H^{-1}t}>x+ct\big)P\big((1-\delta)\mu_H^{-1}t\le N(t)\le(1+\delta)\mu_H^{-1}t\big)\nonumber\\
&\sim&\mu_H^{-1}t\overline{G}\Big((x-\mu_G\mu_H^{-1}t+ct)\big(1+\delta\mu_H^{-1}t(x-\mu_G\mu_H^{-1}t+ct)^{-1}\big)\Big)\nonumber\\
&\ge&\mu_H^{-1}t\overline{G}\Big((x-\mu_G\mu_H^{-1}t+ct)\big(1+\delta\mu_H^{-1}(\gamma-\mu_G\mu_H^{-1}+c)^{-1}\big)\Big)\nonumber\\
&\gtrsim&\mu_H^{-1}t\overline{G}(x-\mu_G\mu_H^{-1}t+ct)\ \ \ \ \ \text{then let}\ \ \delta\downarrow0.
\end{eqnarray}

Further, when $t\to\infty$, for each $\gamma>0$, it holds uniformly for $x\ge\gamma t$ that
\begin{eqnarray}\label{53352}
&&\overline{G}(x-\mu_G\mu_H^{-1}t+ct)=\overline{G}\Big((x-\mu_{G_2}\mu_H^{-1}t+ct)
\big(1+\mu_H^{-1}(\mu_{G_2}-\mu_G)t(x-\mu_{G_2}\mu_H^{-1}t+ct)^{-1}\big)\Big)\nonumber\\
&\ge&\overline{G}\Big((x-\mu_{G_2}\mu_H^{-1}t+ct)
\big(1+\mu_H^{-1}(\mu_{G_2}-\mu_G)(\gamma-\mu_{G_2}\mu_H^{-1}+c)^{-1}\big)\Big)\nonumber\\
&\gtrsim&\overline{G}(x-\mu_{G_2}\mu_H^{-1}t+ct)
\overline{G_*}\big(1+\mu_H^{-1}(\mu_{G_2}-\mu_G)(\gamma-\mu_{G_2}\mu_H^{-1}+c)^{-1}\big).
\end{eqnarray}

Combining (\ref{53352}) and (\ref{53351}), we know that,
for each $\gamma>\max\{1,\mu_G\}$,
\begin{eqnarray}\label{50902}
\liminf\limits_{t\to\infty}\inf\limits_{x\ge\gamma t}
\psi(x;t)\big(\mu_H^{-1}t\overline{G}(x-\mu_{G_2}\mu_H^{-1}t+ct)\big)^{-1}\ge a.
\end{eqnarray}
and
\begin{eqnarray}\label{50903}
\liminf\limits_{t\to\infty}\inf\limits_{x\ge\gamma t}
\psi(x;t)\big(\mu_H^{-1}t\overline{G}(x-\mu_{G}\mu_H^{-1}t+ct)\big)^{-1}\ge1.
\end{eqnarray}

From (\ref{50901}), (\ref{50900}), (\ref{50902}) and (\ref{50903}),
we know that (\ref{509000}) and (\ref{5090000}) hold for each $\gamma>\max\{1,\mu_G\}$.
\\

$(ii)$ Based on the proof of $(i)$, in order to prove (\ref{509001}) for each $\gamma>0$,
we only to get the uniform asymptotic upper bound of $\psi_{2}(x;t)$.
In fact, there is a $t_0>0$ such that $\gamma s(t)>\mu_H^{-1}\mu_G$ for all $t\ge t_0$
and some positive function $s(\cdot)$ on $[0,\infty)$ satisfying $s(t)\to\infty$ as $t\to\infty$.
Further, if $x\ge\gamma ts(t)$, then we have
\begin{eqnarray*}%\label{3308000}
&&\psi_{2}(x;t)\sim(1+\delta)\mu_H^{-1}t\overline{G}\big((x-\mu_H^{-1}\mu_Gt+ct)
(1-\delta-(1+\delta)ct(x-\mu_H^{-1}\mu_Gt+ct)^{-1}\big)\nonumber\\
&\le&(1+\delta)\mu_H^{-1}t\overline{G}\big((x-\mu_H^{-1}\mu_Gt+ct)
(1-\delta-(1+\delta)c(\gamma s(t)-\mu_H^{-1}\mu_G+c)^{-1}\big)\nonumber\\
&\lesssim&(1+\delta)\mu_H^{-1}t\overline{G}(x-\mu_H^{-1}\mu_Gt+ct)\overline{G_*}^{-1}(1-2\delta)\nonumber\\
&\sim&\mu_H^{-1}t\overline{G}(x-\mu_H^{-1}\mu_Gt+ct)\ \ \ \ \ \ \ \ \text{then let}\ \ \delta\downarrow0.
\end{eqnarray*}

\subsection{\normalsize\bf Proof of Theorem \ref{thm1020}}

Based on the proofs of (\ref{3304}) with requirement $\gamma>\mu_G$ and (\ref{3307000}) with requirement
$\gamma>\mu_H^{-1}\mu_G$ in $(i)$,
in order to prove (\ref{5090002}) for each pair $\max\{\mu_G,\mu_H^{-1}\mu_G\}<\gamma<\Gamma<\infty$,
we just need to get the uniformly asymptotic upper bound of $\psi_{3}(x;t)$,
and the uniformly asymptotic lower bound of $\psi(x;t)$ in Case 2.

%For $\psi_{21}(x;t)$, because $(Y_i,Z_i),i\ge1$ are independent random vectors, %according to Lemma \ref{lem3302},
%we know that $Y_i+cZ_i,i\ge1$ still are independent.
%Therefore, according to Lemma \ref{lem3301}, (\ref{3308}) still holds.

For the former, according to Markov inequality and Proposition \ref{pron106}
(thus Proposition \ref{pron105} and Proposition \ref{pron104}),
for each pair $0<\gamma<\Gamma<\infty$ and each $r>J_G^++1$, when $t\to\infty$,
it holds uniformly for all $x\in[\gamma t,\Gamma t]$ that
\begin{eqnarray}\label{331200}
&&\psi_3(x,t)\le P\big(N(t)>(1+\delta)\mu^{-1}_Ht\big)\nonumber\\
&=&O\big(\mu_H^{-1}t^{-r}EN^r(t)\textbf{1}_{\{N(t)>(1+\delta)\mu^{-1}_Ht\}}\big)\nonumber\\
&=&o\big(\mu_H^{-1}t\overline{G}(t)\big)\nonumber\\
&=&o\big(\mu_H^{-1}t\overline{G}(\Gamma^{-1}x)\big)\nonumber\\
&=&o\big(\mu_H^{-1}t\overline{G}(x-\mu_G\mu_H^{-1}t+ct)\big).
\end{eqnarray}

For the latter, according to Lemma \ref{lem201} and Proposition \ref{pron106}, by $\sum_{i=1}^{N(t)}Z_i\le t$, we have
\begin{eqnarray}\label{33141}
&&\psi(x;t)\ge P\Big(\sum_{i=1}^{N(t)}(Y_i-cZ_i)>x\Big)\nonumber\\
&\ge&P\Big(\sum_{i=1}^{N(t)}Y_i>x+ct,(1-\delta)\mu_H^{-1}t\le N(t)\le(1+\delta)\mu_H^{-1}t\Big)\nonumber\\
&\ge&P\big(S^G_{(1-\delta)\mu_H^{-1}t}>x+ct\big)-P\big(|N(t)\mu_Ht^{-1}-1|>\delta\big)\nonumber\\
&=&\psi_{01}(x;t)-\psi_{02}(x;t).
\end{eqnarray}
For $\psi_{01}(x;t)$, by $x\ge\gamma t$ with each $\gamma>\mu_G$ and $c>\mu_H^{-1}\mu_G$, we have
$$x+ct\ge(\gamma+c)t=(\gamma+c)(1-\delta)^{-1}\mu_H(1-\delta)\mu_H^{-1}t=\gamma_1(1-\delta)\mu_H^{-1}t$$
for any $\delta\in(0,1)$ and
$$\gamma_1>(\gamma+\mu_H^{-1}\mu_G)\mu_H=\gamma\mu_H+\mu_G>\mu_G.$$
Then, we further require
$$0<\delta<(\gamma-\mu_H^{-1}\mu_G+c)(\mu_H^{-1}\mu_G)^{-1}.$$
Therefore, according to Lemma \ref{lem201} and by $G\in\mathcal{C}$,
for each pari $\mu_G<\gamma<\Gamma<\infty$, when $t\to\infty$, it holds uniformly and for all $x\in[\gamma t,\Gamma t]$ that
\begin{eqnarray}\label{33151}
&&\psi_{01}(x;t)\sim(1-\delta)\mu_H^{-1}t\overline{G}\big(x-(1-\delta)\mu_H^{-1}\mu_Gt+ct\big)\nonumber\\
&=&(1-\delta)\mu_H^{-1}t\overline{G}\big((x-\mu_H^{-1}t\mu_G+ct)(1+\delta\mu_H^{-1}\mu_Gt(x-\mu_H^{-1}t\mu_G+ct)^{-1}\big)\nonumber\\
&\ge&(1-\delta)\mu_H^{-1}t\overline{G}\big((x-\mu_H^{-1}t\mu_G+ct)(1+\delta\mu_H^{-1}\mu_G(\Gamma-\mu_H^{-1}\mu_G+c)^{-1}\big)\nonumber\\
&\sim&\mu_H^{-1}t\overline{G}\big((x-\mu_H^{-1}t\mu_G+ct)\big)\ \ \ \ \ \ \ \ \ \text{then let}\ \ \ \ \delta\downarrow0.
\end{eqnarray}
For $\psi_{02}(x;t)$, by (4.7) of Chen et al. \cite{CWY2021}, it holds uniformly and for all $x\in[\gamma t,\Gamma t]$ that
\begin{eqnarray}\label{33161}
\psi_{02}(x;t)=o\big(\mu_H^{-1}t\overline{G}(x)\big)=o\big(\mu_H^{-1}t\overline{G}\big(x-\mu_H^{-1}t\mu_G+ct)\big).
\end{eqnarray}
Combining (\ref{331200})-(\ref{33161}), we can prove that, for each pair $\mu_G<\gamma<\Gamma<\infty$,
when $t\to\infty$, it holds uniformly for all $x\in[\gamma t,\Gamma t]$ that
\begin{eqnarray*}  %\label{33181}
&&\psi(x;t)\gtrsim\mu_H^{-1}t\overline{G}\big(x-\mu_H^{-1}t\mu_G+ct).
\end{eqnarray*}

In this way, we have completed the proof of the theorem.

\subsection{\normalsize\bf Proof of Theorem \ref{thm102}}

$(i)$ Using any integer $k$ and any $\delta$ such as the proof of Theorem \ref{thm503} $(i)$, we split
\begin{eqnarray}\label{3317}
P\big(R_{0}(t)>x\big)&=&\bigg(\sum_{m=1}^{\lfloor(1+\delta)t(k\mu_H)^{-1}\rfloor}
+\sum_{m=\lfloor(1+\delta)t(k\mu_H)^{-1}\rfloor+1}^{\lfloor(1+\delta)t\mu_H^{-1}\rfloor}
+\sum_{m=\lfloor(1+\delta)t\mu_H^{-1}\rfloor+1}^{\infty}\bigg)P\big(S_m^F>x,N(t)=m\big)\nonumber\\
&=&P_1(x;t)+P_2(x;t)+P_3(x;t).
\end{eqnarray}
By the proof of Theorem \ref{thm503} $(i)$, we just need to get the uniform asymptotic upper bound of $P_2(x;t)$
for $x\ge\gamma t$ with each $\gamma>\mu_G$.
To this end, considering $S^H_{m+1}>t$ when $N(t)=m$, we further split
\begin{eqnarray}\label{3318}
&&P_2(x;t)\le\sum_{\lfloor(1-\delta)\mu^{-1}_Ht\rfloor\le m\le\lfloor(1+\delta)\mu^{-1}_Ht\rfloor}
P\big(S^G_m-cS^H_{m+1}+cZ_{m+1}>x,N(t)=m\big)\nonumber\\
&\le&P\big(S^G_{\lfloor(1+\delta)\mu^{-1}_Ht\rfloor}+\max\{cZ_i,1\le i\le \lfloor(1+\delta)\mu^{-1}_Ht\rfloor+1\}>x+ct\big)\nonumber\\
&\le&P\big(S^G_{\lfloor(1+\delta)\mu^{-1}_Ht\rfloor}>(1-\delta)(x+ct)\big)
+P\big(\max\{cZ_i,1\le i\le \lfloor(1+\delta)\mu^{-1}_Ht\rfloor+1\}>\delta(x+ct)\big)\nonumber\\
&=&P_{21}(x;t)+P_{22}(x;t).
\end{eqnarray}

We first deal with $P_{21}(x;t)$. Further set
$$0<\delta<\min\{\gamma\mu_H(\gamma\mu_H+2\mu_G)^{-1},(\gamma-\mu_G\mu_H^{-1}+c)(\gamma+\mu_G\mu_H^{-1}+c)^{-1}\},$$
by $\gamma>\mu_G$ and $x\ge\gamma t$, we have
$$(1-\delta)(x+ct)\ge(1-\delta)(\gamma+c)t=(1-\delta)(1+\delta)^{-1}(\gamma+c)\mu_H(1+\delta)\mu^{-1}_Ht=\gamma_1(1+\delta)\mu^{-1}_Ht$$
and
$$\gamma_1>(1-\delta)(1+\delta)^{-1}(\gamma+\mu^{-1}_H\mu_G)\mu_H=(1-\delta)(1+\delta)^{-1}(\gamma\mu_H+\mu_G)>\mu_G.$$
Then when $t\to\infty$, according to Lemma \ref{lem201}, by $G\in\mathcal{C}$,
we know that, for each $\gamma>\mu_G$, it holds uniformly and for all $x\ge\gamma t$ that
\begin{eqnarray}\label{3319}
&&P_{21}(x;t)\lesssim(1+\delta)\mu^{-1}_Ht\overline{G}\big((1-\delta)(x+ct)-(1+\delta)\mu_H^{-1}\mu_Gt\big)\nonumber\\
&=&(1+\delta)\mu^{-1}_Ht\overline{G}\big((x-\mu_G\mu_H^{-1}t+ct)-\delta(x+\mu_H^{-1}\mu_Gt+ct)\big)\nonumber\\
&=&(1+\delta)\mu^{-1}_Ht\overline{G}\big((x-\mu_G\mu_H^{-1}t+ct)(1-\delta-2\delta\mu_H^{-1}\mu_Gt(x-\mu_G\mu_H^{-1}t+ct)^{-1})\big)\nonumber\\
&\le&(1+\delta)\mu^{-1}_Ht\overline{G}\big((x-\mu_G\mu_H^{-1}t+ct)(1-\delta-2\delta\mu_H^{-1}\mu_G(\gamma-\mu_G\mu_H^{-1}+c)^{-1})\big)\nonumber\\
&\sim&\mu_H^{-1}t\overline{G}(x-\mu_G\mu_H^{-1}t+ct)\ \ \ \ \ \ \ \ \text{then let}\ \ \ \delta\downarrow0.
\end{eqnarray}

Next, we deal with $P_{22}(x;t)$.
By (\ref{532}) and $G\in\mathcal{C}$,
we know that, when $t\to\infty$, for each $\gamma>\mu_G$, it is holds uniformly for all $x\ge\gamma t$ that
\begin{eqnarray}\label{3320}
&&P_{22}(x;t)\le(1+\delta)(\mu^{-1}_Ht+1)\overline{H}(c^{-1}\delta(x+ct)\big)\nonumber\\
&\le&(1+\delta)(\mu^{-1}_Ht+1)\overline{H}(c^{-1}\delta(x-\mu_G\mu_H^{-1}+ct)\big)\nonumber\\
&=&o\big(\mu^{-1}_Ht\overline{G}(x-\mu_G\mu_H^{-1}+ct)\big).
\end{eqnarray}

Combined with (\ref{3317})-(\ref{3320}), (\ref{10015}) holds for each $\gamma>\mu_G$.
\\

$(ii)$ To prove (\ref{203}),
the proof of the uniform asymptotic upper bounds for $R_0(t)$, see the proof of $(i)$ of the theorem.
And the proof of the uniform asymptotic lower bounds for $R_0(t)$, see the proof of Theorem \ref{thm1020}.

\subsection{\normalsize\bf Proof of Theorem \ref{thm602}}

Clearly, by $\mu_G<\infty$ and (\ref{408}), we have $E\tau<\infty$.
In order to prove $EN(\tau)<\infty$, we need the following lemma.

\begin{lemma}\label{lem302}
In the above nonstandard renewal risk model,
let $Z_i,\ i\ge1$ be WLOD random variables with dominating coefficients $g_{L,H}(n),n\ge1$ satisfying (\ref{105}).
Then $E\tau<\infty$ if and only if $EN(\tau)<\infty$.
\end{lemma}

\proof According to Proposition \ref{pron105}, by (\ref{105}),
we know that for any $\delta\in(0,1)$, there exist a $t_0=t_0(H,\delta)>0$ such that,
for all $t\ge t_0$, it holds that
\begin{eqnarray}\label{3012}
(1-\delta)\mu^{-1}_Ht<EN(t)<(1+\delta)\mu^{-1}_Ht.
\end{eqnarray}
If $E\tau<\infty$, then by (\ref{3012}), we have
\begin{eqnarray}\label{3013}
&&EN(\tau)=\Big(\int_{0-}^{t_0}+\int_{t_0}^\infty\Big)EN(t)P(\tau\in dt)\nonumber\\
&\le&\int_{0-}^{t_0}EN(t_0)P(\tau\in dt)+\int_{t_0}^\infty EN(t)P(\tau\in dt)\nonumber\\
&<&(1+\delta)\mu^{-1}_H\int_{0-}^\infty(t_0+t)P(\tau\in dt)<\infty.
\end{eqnarray}
Similarly, if $EN(\tau)<\infty$, we can prove $E\tau<\infty$.
\hfill$\Box$

According to Lemma \ref{lem302}, by $E\tau<\infty$, we know that $EN(\tau)<\infty$.
\\

Now, we prove (\ref{409}). To this end, we take any $\gamma>\max\{\mu_G,\mu_H^{-1}\mu_G\}$ to split
\begin{eqnarray}\label{311}
&&\psi(x;\tau)=\Big(\int_{0-}^{\gamma^{-1}x}+\int_{\gamma^{-1}x}^\infty\Big)P\big(\max_{1\le k\le N(t)}S^F_k>x\big)P(\tau\in dt)\nonumber\\
&=&I_1(x)+I_2(x),\ \ \ \ \ \ \ \ \ x\in(0,\infty).
\end{eqnarray}

We first deal with $I_2(x)$.  By (\ref{408}), we know that
\begin{eqnarray}\label{313}
I_2(x)\le P(\tau>\gamma^{-1}x)=o\big(\overline{G}(x)\big).
\end{eqnarray}
Therefore, we only need to with $I_1(x)$.

For any $\varepsilon\in(0,1)$ and any $\gamma>\mu_G$,
according to Lemma \ref{lem201},
there exists an integer $n_0=n_0(\varepsilon,G)\ge1$ such that, when $n> n_0$, it holds that
\begin{eqnarray}\label{314}
(1-\varepsilon)\overline{G}_*\Big(1+c\gamma^{-1}\Big)
\le\inf_{x\ge\max\{\gamma t,\gamma n\}}\frac{P(S^G_n>x+ct)}{n\overline{G}(x)}\le\sup_{x\ge\gamma n }\frac{P(S^G_n>x)}{n\overline{G}(x)}
\le\frac{(1+\varepsilon)}{\overline{G_*}\big(\frac{1}{1-\frac{\mu_G}{\gamma})}\big)}.
\end{eqnarray}
%According to Lemma \ref{lem201} again,
By $G\in\mathcal{C}$, for the above $\varepsilon$, $n_0$ and all $n\leq n_0$,
there exists $x_0=x_0(n_0,\varepsilon)>0$ large enough such that
\begin{eqnarray}\label{315}
(1-\varepsilon)\overline{G_*}(1+c\gamma^{-1})\le\inf_{x\ge \max\{\gamma t,x_0\}}\frac{P\big(S^G_n>x+ct\big)}{n\overline{G}(x)}
\le\sup_{x\ge x_0}\frac{P\big(S^G_n>x\big)}{n\overline{G}(x)}\le1+\varepsilon
\end{eqnarray}
In addition, by $EN(\tau)<\infty$,
\begin{eqnarray}\label{40000}
EN(\tau)\textbf{1}_{\{N(\tau)\ge\min\{N(x),\ \gamma^{-1}x\}\}}\to0.
\end{eqnarray}

Then by the above $n_0$, we again split
\begin{eqnarray}\label{316}
I_1(x)=\int_{0-}^{\gamma^{-1}x}\Big(\sum_{n=1}^{n_0}+\sum_{n=n_0+1}^\infty\Big)
P\big(\max_{1\le k\le n}S^F_k>x,\ N(t)=n\big)P(\tau\in dt)=I_{11}(x)+I_{12}(x).
\end{eqnarray}

For $I_{11}(x)$, by (\ref{315}), we have
\begin{eqnarray}\label{317}
&&I_{11}(x)\le\int_{0-}^{\gamma^{-1}x}\sum_{n=1}^{n_0}P(S^G_n>x)P\big(N(t)=n\big)P(\tau\in dt)\nonumber\\
&<&(1+\varepsilon)\int_{0-}^{\gamma^{-1}x}\sum_{n=1}^{n_0}nP\big(N(t)=n\big)P(\tau\in dt)\overline{G}(x)\nonumber\\
&\le&(1+\varepsilon)EN(\tau)\textbf{1}_{\{N(\tau)\le n_0\}}\overline{G}(x)\nonumber\\
\end{eqnarray}
and
\begin{eqnarray}\label{317-2}
&&I_{11}(x)\geq\int_{0-}^{\gamma^{-1}x}\sum_{n=1}^{n_0}P\big(S_n^F>x,N(t)=n\big)P(\tau\in dt)\nonumber\\
&\ge&\int_{0-}^{\gamma^{-1}x}\sum_{n=1}^{n_0}P\big(S_n^G>x+ct)P\big(N(t)=n\big)P(\tau\in dt)\nonumber\\
%&>&(1-\varepsilon)^{2^{-1}}EN(\tau)\textbf{1}_{\{N(\tau)\le n_0\}}\overline{G}\big(x(1+ctx^{-1})\big)\nonumber\\
&>&(1-\varepsilon)EN(\tau)\textbf{1}_{\{N(\tau)\le n_0\}}\overline{G_*}(1+c\gamma^{-1})\overline{G}(x)-P(\tau>\gamma^{-1}x).
\end{eqnarray}

For $I_{12}(x)$, we further split it as follows:
\begin{eqnarray}\label{318}
&&I_{12}(x)=\int_{0-}^\infty\sum_{n=n_0+1}^\infty\big(\textbf{1}_{\{t:{\gamma^{-1}x}\ge\max\{t,n\}}\}
+\textbf{1}_{\{t:0\le t\le{\gamma^{-1}x}<n}\}\big)\nonumber\\
&&\ \ \ \ \ \ \ \ \ \ \ \ \ \ \ P\big(\max_{1\le k\le n}S_k^F>x,N(t)=n\big)P(\tau\in dt)\nonumber\\
&=&I_{121}(x)+I_{122}(x).
\end{eqnarray}
By (\ref{314}) and $c\mu_H>\mu_G$, we have
\begin{eqnarray}\label{319}
&&I_{121}(x)\le\int_{0-}^\infty\sum_{n=n_0+1}^\infty P(S_n^G>x)P\big(N(t)=n\big)
\textbf{1}_{\{t:\gamma^{-1}x\ge\max\{t, n\}\}}(t)P(\tau\in dt)\nonumber\\
%&\lesssim&\int_{0-}^\infty\sum_{n=n_0+1}^\infty n\overline{G}(x-n\mu_G)P\big(N(t)=n\big)
%\textbf{1}_{\{t:\gamma^{-1}x\ge\max\{t,\mu_Hn\}\}}(t)P(\tau\in dt)\nonumber\\
%&\le&\int_{0-}^\infty\sum_{n=n_0+1}^\infty n\overline{G}\Big(x\big(1+(c-\mu_H^{-1}\mu_G)\gamma^{-1}\big)\Big)P\big(N(t)=n\big)
%\textbf{1}_{\{t:\gamma^{-1}x\ge\max\{t,\mu_Hn\}\}}(t)P(\tau\in dt)\nonumber\\
&\le&\frac{(1+\varepsilon)\overline{G}(x)}{\overline{G_*}\big((1-\mu_G\gamma^{-1})^{-1}\big)}
\bigg(\int_{0-}^\infty\sum_{n=n_0+1}^\infty nP\big(N(t)=n\big)P(\tau\in dt)\nonumber\\
&&\ \ \ \ \ \ \ \ \ \ \ \ \ \ \
-\int_{0-}^\infty\sum_{n=n_0+1}^\infty nP\big(N(t)=n\big)\textbf{1}_{\{t:0\le t\le\gamma^{-1}x< n}\}(t)P(\tau\in dt)\bigg)\nonumber\\
&=&(1+\varepsilon)\overline{G_*}^{-1}\big((1-\mu_G\gamma^{-1})^{-1}\big)
\big(EN(\tau)\textbf{1}_{\{N(\tau)>n_0\}}-I_{120}(x)\big)\overline{G}(x),
\end{eqnarray}
\begin{eqnarray}\label{3190}
&&I_{121}(x)\ge\int_{0-}^\infty\sum_{n=n_0+1}^\infty\textbf{1}_{\{t:\gamma^{-1}x\ge\max\{t, n\}\}}
P(S_n^G>x+ct)P\big(N(t)=n\big)P(\tau\in dt)\nonumber\\
&\ge&(1-\varepsilon)\overline{G_*}(1+c\gamma^{-1})\overline{G}(x)\int_{0-}^\infty\sum_{n=n_0+1}^\infty nP\big(N(t)=n\big)\textbf{1}_{\{t:\gamma^{-1}x\ge\max\{t, n\}\}}(t)P(\tau\in dt)\nonumber\\
&=&(1-\varepsilon)\overline{G_*}(1+c\gamma^{-1})\overline{G}(x)
\big(EN(\tau)\textbf{1}_{\{N(\tau)>n_0,x\ge\gamma\tau\}}-I_{120}(x)\big)\nonumber\\
&\ge&(1-\varepsilon)\overline{G_*}(1+c\gamma^{-1})\overline{G}(x)
\big(EN(\tau)\textbf{1}_{\{N(\tau)>n_0\}}-EN(\tau)\textbf{1}_{\{\tau>\gamma^{-1}x\}}-I_{120}(x)\big)
\end{eqnarray}
and
\begin{eqnarray}\label{435}
&&I_{122}(x)\le(1+\varepsilon)\overline{G}^{-1}_*\big((1-\mu_G\gamma^{-1})^{-1}\big)\overline{G}(x)I_{120}(x).
\end{eqnarray}
Here, by (\ref{40000}),
\begin{eqnarray}\label{320}
I_{120}(x)\le EN(\tau)\textbf{1}_{\{N(\tau)>\gamma^{-1}x}\}\to0.
\end{eqnarray}

Combining (\ref{311}), (\ref{313}) and (\ref{316})-(\ref{320}),
writing
$$b_1(\gamma)=\overline{G_*}(1+c\gamma^{-1})\ \ \ \text{and}\ \ \
b_2(\gamma)=\overline{G}^{-1}_*\Big(\big(1+(c-\mu_H^{-1}\mu_G)\gamma^{-1}\big)^{-1}\Big),$$
%\textcolor[rgb]{0.00,0.07,1.00}{and by (\ref{408})},
and by (\ref{408}), we have
\begin{eqnarray}\label{4120001}
(1-\varepsilon)b_1(\gamma)\le\liminf\frac{\psi(x;\tau)}{EN(\tau)\overline{G}(x)}
\le\limsup\frac{\psi(x;\tau)}{EN(\tau)\overline{G}(x)}\le(1+\varepsilon)b_2(\gamma).
\end{eqnarray}
Let $\gamma\uparrow\infty$ and $\varepsilon\downarrow0$ in (\ref{4120001}), then by $G\in\mathcal{C}$, we know that (\ref{409}) holds.
$\hspace{\fill}\Box$\\

%At the end of this subsection, we give a sufficient condition for (\ref{408}).
%
%\begin(1+\varepsilon){pron}\label{pron107}
%In the above nonstandard renewal risk model Case $(i)$,
%let $Z_i,i\ge1$ be WLOD random variables with dominating coefficients $g_{L,H}(n),\ n\ge1$ satisfying all conditions in Proposition \ref{pron106}.
%If $G\in\mathcal{D}$ and $EN^p(\tau)<\infty$ for some $p>J_G^+$,
%then (\ref{408}) holds.
%\end{pron}
%
%\proof According to Proposition \ref{Proposition102} $(iii)$,
%for each $p>J_G^+$, by $G\in\mathcal{D}$, we have
%$$\lim x^p\overline{G}(x)=\infty.$$
%And by $EN^p(\tau)<\infty$, we have
%$$\lim x^pP\big(N(\tau)>x\big)=0.$$
%Combined with the above two asymptotic expressions, we know that
%$$P\big(N(\tau)>x\big)=o\big(\overline{G}(x)\big).$$
%By the above fact, we have
%\begin{eqnarray*}
%&&P(\tau>x)=P\big(\tau>x,N(\tau)>2^{-1}\mu_H^{-1}x\big)+P\big(\tau>x,N(\tau)\le2^{-1}\mu_H^{-1}x\big)\\
%&\le&P\big(N(\tau)>2^{-1}\mu_H^{-1}x\big)+\int_{x}^\infty P\big(N(t)\le2^{-1}\mu_H^{-1}x\big)P(\tau\in dt)\\
%&\le&o\big(\overline{G}(x)\big)+P\big(N(x)\le2^{-1}\mu_H^{-1}x\big)P(\tau>x).
%\end{eqnarray*}
%Therefore, for $x$ large enough, it holds that
%$$P(\tau>x)\le o\Big(\overline{G}(x)P^{-1}\big(N(x)>2^{-1}\mu_H^{-1}x\big)\Big)=o\big(\overline{G}(x)\big)$$
%according to Proposition \ref{pron106}.
%$\hspace{\fill}\Box$

\section{\normalsize\bf Some applications}
\setcounter{equation}{0}\setcounter{thm}{0}\setcounter{pron}{0}\setcounter{lemma}{0}\setcounter{Corol}{0}

In this section, we utilize Theorem \ref{thm102} to investigate the precise deviation of the proportional net loss process
and excess-of-net-loss process, as well as the asymptotic estimate of the mean of stop-net-loss reinsurance treaty.
These studies are all crucial components in risk theory.
For the sake of brevity, we have omitted the direct proof of these results.%, which can be obtained directly from Theorem \ref{thm102}.}

\subsection{\normalsize\bf Precise large deviations of the proportional net loss process}

In a renewal risk model, at time $t\ge0$,
a reinsurance company covers percent $q_1\cdot 100$ of the total claim amount $S_{N(t)}^G=\sum_{i=1}^{N(t)}Y_i$ for some $q_1\in(0,1]$;
the reinsurance company benefits percent $q_2\cdot 100$ of the total income $cS_{N(t)}^H=c\sum_{i=1}^{N(t)}Z_i$ for some $q_2\in[0,1]$.
Then we respectively introduce the concepts of proportional reinsurance process
and proportional net loss process for the reinsurance company defined by
$$\big(R_{001}(t)=q_1S_{N(t)}^G:t\ge0\big)$$
and
$$\Big(R_{01}(t)=q_1S_{N(t)}^G-q_2cS_{N(t)}^H=q_1\sum_{i=1}^{N(t)}\big(Y_i-q_1^{-1}q_2 cZ_i\big)
=q_1\sum_{i=1}^{N(t)}X_i=q_1S_{N(t)}^{F_1}:t\ge0\Big),$$
where $F_1$ is the distribution of $X_1$ and $S_{n}^{F_1}=\sum_{i=1}^nX_i,\ n\ge1$.
Clearly, $F_1$ is supported on $(-\infty,\infty)$.
Naturally, we require the following more general safety load condition to be satisfied:
$$q_1\mu_G<q_2c\mu_H.$$
In particular, when $q_1=q_2=1$, $F_1=F$ and $q_1\mu_G<q_2c\mu_H$ is $\mu_G<c\mu_H$.

The research goal $R_{001}(t)$ comes from Example 5.2 $(i)$ of Kl\"{u}ppelberg and Mikosch $\cite{KM1997}$,
which have obtained the following conclusion. For each $\gamma>q_1\mu_G$,
\begin{eqnarray*}
\lim\limits_{t\to\infty}\sup\limits_{x\ge\gamma\lambda(t)}
\Big|P\big(R_{001}(t)>x\big)\Big(\lambda(t)\overline{G}\big(q_1^{-1}x-\mu_G\lambda(t)\big)\Big)^{-1}-1\Big|=0
\end{eqnarray*}
under the assumption that %$Y_i,\ i\ge1$ are independent random variables with common distribution
$G\in\mathcal{ERV}(\alpha,\beta)$ for some pair $1<\alpha\le\beta<\infty$ in a standard renewal risk model.

In this subsection, we focus on $R_{01}(t)$ in a nonstandard renewal risk model.
Clearly, $R_{01}(t)$ is naturally inspired by $R_{001}(t)$.
However, the former has a wider application, see Theorem \ref{thm101} below.

Especially, when $q_1=q_2=1$, $R_{001}(t)=R_{00}(t)$ and $R_{01}(t)=R_{0}(t),\ t\ge0$.

\begin{thm}\label{thm101}
$(i)$ In Case 1, conditions such as Theorem \ref{thm102} $(i)$, then for each $\gamma>q_1\mu_G$,
\begin{eqnarray}\label{426}
\lim\limits_{t\to\infty}\sup\limits_{x\ge\gamma t}\Big|P\big(R_{01}(t)>x\big)
\Big(\mu_H^{-1}t\overline{G}\big(q_1^{-1}x-\mu_G\mu_H^{-1}t+q_2cq_1^{-1}t\big)\Big)^{-1}-1\Big|=0.
\end{eqnarray}

$(ii)$ In Case 2, conditions such as Theorem \ref{thm102} $(ii)$, then for each pair $q_1\mu_G<\gamma<\Gamma<\infty$,
\begin{eqnarray}\label{427}
\lim\limits_{t\to\infty}\sup\limits_{x\in[\gamma t,\Gamma t]}\Big|P\big(R_{01}(t)>x\big)
\Big(\mu_H^{-1}t\overline{G}\big(q_1^{-1}x-\mu_G\mu_H^{-1}t+q_2cq_1^{-1}t\big)\Big)^{-1}-1\Big|=0.
\end{eqnarray}
\end{thm}

We omit the proof details of this theorem, which is similar to that of Theorem \ref{thm102}.

\subsection{\normalsize\bf Precise large deviations of excess-of-net-loss process}

We say the stochastic process
$$\Big(R_{002}(t)=\sum\limits_{i=1}^{N(t)}(Y_i-D)^+:t\ge0\Big)$$
an excess-of-loss process, where $D$ is some positive finite constant,
see Example 5.2 $(ii)$ of Kl\"{u}ppelberg and Mikosch $\cite{KM1997}$.
Additionally, we believe that introducing income factors will be better to evaluating reinsurance
risks and determining more reasonable threshold $D$.
Then, we focus on the corresponding research objects excess-of-net-loss process
$$\Big(R_{02}(t)=\sum_{i=1}^{N(t)}\big((Y_i-D)^+-q_2 cZ_i\big):t\ge0\Big).$$

To give three conclusions for $R_{22}(t)$, let $G_1$ be distribution of $(Y_1-D)^+$. Thus
$$\overline{G_1}(x)=\textbf{1}_{(-\infty,0)}(x)+\overline{G}(x+D)\textbf{1}_{[0,\infty)}(x),\ \ \ x\in(-\infty,\infty)$$
and
$$\mu_{G_1}=E(Y_1-D)^+=EY_1\textbf{1}_{\{Y_1>D\}}-DP(Y_1>D).$$
\\

\begin{thm}\label{cor2001}
$(i)$ In Case 1, under conditions as Theorem \ref{thm102} $(i)$, for each $\gamma>\mu_{G_1}$,
\begin{eqnarray}\label{204}
\lim\limits_{t\to\infty}\sup\limits_{x\ge\gamma t}\big|P\big(R_{02}(t)>x\big)
\big(\mu_H^{-1}t\overline{G}(x+D-\mu_{G_1}\mu_H^{-1}t+q_2ct)\big)^{-1}-1\big|=0.
\end{eqnarray}

$(ii)$ In Case 2, under conditions as Theorem \ref{thm102} $(ii)$, for each pair $\mu_{G_1}<\gamma<\Gamma<\infty$,
\begin{eqnarray}\label{205}
\lim\limits_{t\to\infty}\sup\limits_{x\in[\gamma t,\Gamma t]}\big|P\big(R_{02}(t)>x\big)
\big(\mu_H^{-1}t\overline{G}(x+D-\mu_{G_1}\mu_H^{-1}t+q_2ct)\big)^{-1}-1\big|=0.
\end{eqnarray}
\end{thm}

The proof is also similar to that of Theorem \ref{thm102}.

\subsection{\normalsize\bf Asymptotic estimate of mean of stop-loss reinsurance process}

We denote stop-loss reinsurance process
$$\bigg(R_{003}(t)=\Big(\sum_{i=1}^{N(t)}Y_i\big(c\mu_G\lambda(t)\big)^{-1}-K\Big)^+:\ t\ge0\bigg),$$
where $K>c^{-1}$ and $c$ is the premium income mentioned above,
see Example 5.2 $(iii)$ of Kl\"{u}ppelberg and Mikosch $\cite{KM1997}$.
However, we believe that the premium income $c$ has no obvious significance here.
Therefore, we may as well remove this factor and replace $\lambda(t)$ with $\mu_H^{-1}t$ in the following.
Further, for $q_1\in(0,1]$ and $q_2\in[0,1]$, by the safety load condition, we have
$$\mu_0=|q_2\mu_H-q_1c\mu_G|>0.$$
Then we focus on the corresponding research objects stop-net-loss reinsurance process
\begin{eqnarray}\label{207}
\bigg(R_{03}(t)=\Big(\sum_{i=1}^{N(t)}(q_1Y_i-q_2cZ_i)(\mu_0\mu_H^{-1}t)^{-1}-K\Big)^+:\ t\ge0\bigg).
\end{eqnarray}
To this end, we recall the integral tail distribution of distribution $G$
$$G_I(x)=\mu_G^{-1}\int_{0-}^x\overline{G}(y)dy\textbf{1}_{[0,\infty)}(x),\ \ \ \ x\in(-\infty,\infty),$$
see Kl\"{u}ppelberg $\cite{K1988}$.

\begin{thm}\label{cor2002}
In Case 1 of the above nonstandard renewal risk model, assume that %$K>c^{-1}q_1\mu_0^{-1}\mu_G$
when $q_2=0$ or $K>q_1\mu_0^{-1}\mu_H\mu_G$.
Under conditions as Theorem \ref{thm102} $(i)$, it holds that
\begin{eqnarray}\label{208}
\lim\limits_{t\to\infty}\frac{ER_{03}(t)}{q_1\mu_0^{-1}\mu_G\overline{G_I}(q_1^{-1}\mu_0\mu_H^{-1}Kt-\mu_G\mu_H^{-1}t+q_1^{-1}q_2ct)}=1.
\end{eqnarray}
\end{thm}

\begin{remark}\label{rem203}
$(i)$ In particular, if $q_1=1$ and $q_2=0$,
then $\mu_0=\mu_G$ and $R_{03}(t)=R_{003}(t)$.
Further, if $G\in\mathcal{R}_\alpha$ for some $\alpha>1$,
then according to Theorem \ref{cor2002} and Karamata Theorem,
we can respectively obtain the corresponding results in Corollary 3.4 of Kl\"{u}ppelberg and Mikosch \cite{KM1997}
and Corollary 2.1 of Chen et al. \cite{CWY2021} in more general Case 2.

$(ii)$ Theorem \ref{cor2002} also contains some new results than Corollary 3.4 of Kl\"{u}ppelberg and Mikosch \cite{KM1997}
and Corollary 2.1 of Chen et al. \cite{CWY2021}, in which $\alpha>1$ is required if $G\in\mathcal{R}_{\alpha}$.
\begin{exam}\label{exam2001}
In Theorem \ref{cor2002}, let %$G$ be a distribution of $Y_1$ such that
$$\overline{G}(x)=\emph{\textbf{1}}_{(-\infty,0)}(x)+(1+x)^{-1}\ln^{-2}(e+x)\emph{\textbf{1}}_{[0,\infty)}(x),\ \ \ x\in(-\infty,\infty).$$
Clearly, $G\in\mathcal{R}_1$, $\mu_G\in(0,\infty)$, $EY_1^\alpha=\infty$ for each $\alpha>1$ and
$$\overline{G}(x)\sim x^{-1}\ln^{-2}x.$$
Then according to Theorem \ref{cor2002}, when $t\to\infty$, we have
\begin{eqnarray*}
ER_{03}(t)\sim\frac{q_1}{\mu_0}\int_{\frac{\mu_0K}{q_1\mu_H}-\frac{\mu_Gt}{\mu_H}+\frac{q_2ct}{q_1}}^\infty \frac{1}{v\ln^2v}dv
=\frac{q_1\mu_0^{-1}}{\ln\big((\frac{\mu_0K}{q_1\mu_H}-\frac{\mu_G}{\mu_H}+\frac{q_2c}{q_1})t\big)}
\sim\frac{q_1}{\mu_0}\overline{G}(t)t\ln t.%\asymp ER_{003}(t).
\end{eqnarray*}
\end{exam}
Here, we find an interesting phenomenon that $c,\ K$ and $H$
do not affect the asymptotic expression of $R_{003}(t)$ and $R_{03}(t)$.
For the case that $\alpha>1$, however, they all play a role in the above expression.
%This may be caused by $\alpha=1$ and the properties of slowly varying function $\ln(\cdot)$.
In addition, when $\alpha=1$, $R_{003}(t),\ R_{03}(t)$ and $\overline{G}(t)t\ln t$ are infinitesimals of the same order,
which is different from the case $\alpha\in(1,\infty)$ with the same order as $\overline{G}(t)t$, as $t\to\infty$.
%see Corollary 3.4 of Kl\"{u}ppelberg and Mikosch \cite{KM1997} and Corollary 2.1 of Chen et al. \cite{CWY2021}.
\end{remark}

%The proof of Theorem \ref{cor2002} is similar to that of Corollary 3.4 of Kl\"{u}ppelberg and Mikosch \cite{KM1997},
%obtained by applying Theorem \ref{thm102} of this paper.

%\proof We only prove $(i)$.
%By the safety load condition that
%$q_1\mu_G<q_2c\mu_H$ and $K>q_1\mu_0^{-1}\mu_G$ when $q_2=0$,
%$$\lim_{t\to\infty}q_1^{-1}\mu_0\mu_H^{-1}Kt-\mu_G\mu_H^{-1}t+q_1^{-1}q_2ct=\infty.$$
%Then according to Theorem \ref{thm101}, we have
%\begin{eqnarray*}\label{422}
%&&ER_{03}(t)=\int_{0-}^\infty P\bigg(\sum_{i=1}^{N(t)}(q_1Y_i-q_2cZ_i)>\mu_0\mu_H^{-1}t(K+u)\bigg)du\nonumber\\
%&\sim&\mu^{-1}_H t\int_{0-}^\infty\overline{G}\big(\mu_0\mu_H^{-1}tq_1^{-1}(K+u)-\mu_G\mu_H^{-1}t+q_1^{-1}q_2ct\big)du\nonumber\\
%&=&(q_1^{-1}\mu_0)^{-1}\int_{q_1^{-1}\mu_0\mu_H^{-1}Kt-\mu_G\mu_H^{-1}t+q_1^{-1}q_2ct}^\infty\overline{G}(v)dv\\
%&=&q_1\mu_0^{-1}\mu_G\mu_G^{-1}\int_{q_1^{-1}\mu_0\mu_H^{-1}Kt-\mu_G\mu_H^{-1}t+q_1^{-1}q_2ct}^\infty\overline{G}(v)dv,\ \ \ \ \ \text{as}\ \ t\to\infty\\
%&=&q_1\mu_0^{-1}\mu_G\overline{G_I}(q_1^{-1}\mu_0\mu_H^{-1}Kt-\mu_G\mu_H^{-1}t+q_1^{-1}q_2ct),
%\end{eqnarray*}
%that is (\ref{208}) is holds.
%$\hspace{\fill}\Box$

\end{document}